\documentclass[11pt,reqno]{amsart}
\usepackage[margin=1in]{geometry}
\makeatletter
\def\section{\@startsection{section}{1}%
	\z@{.7\linespacing\@plus\linespacing}{.5\linespacing}%
	{\bfseries
		\centering
}}
\def\@secnumfont{\bfseries}
\makeatother

\usepackage{amsmath,amssymb,amsthm,graphicx,amsxtra, setspace}
\usepackage[utf8]{inputenc}
\usepackage{charter}
\usepackage{mathrsfs}
\usepackage{alltt}
\usepackage{stackengine}
\usepackage{relsize}
\usepackage{hyperref}
\usepackage{aliascnt}
\usepackage{tikz}
\usepackage{mathtools}
\usepackage{multicol}
\usepackage{upgreek}
\usepackage{graphicx,type1cm,eso-pic,color}
\allowdisplaybreaks
\usepackage{scalerel,stackengine}
\stackMath
\usepackage{setspace}
\newtheorem{theorem}{Theorem}[section]
\newtheorem*{theorem*}{Theorem}

\newaliascnt{lemma}{theorem}
\newtheorem{lemma}[lemma]{Lemma}
\aliascntresetthe{lemma} 

\newaliascnt{proposition}{theorem}
\newtheorem{proposition}[proposition]{Proposition}
\aliascntresetthe{proposition}

\newaliascnt{assumption}{theorem}

\aliascntresetthe{assumption}

\newaliascnt{auxiliary}{theorem}

\aliascntresetthe{auxiliary}

\newaliascnt{corollary}{theorem}
\newtheorem{corollary}[corollary]{Corollary}
\aliascntresetthe{corollary}

\newaliascnt{definition}{theorem}
\newtheorem{definition}[definition]{Definition}
\aliascntresetthe{definition}

\newaliascnt{example}{theorem}

\aliascntresetthe{example}

\newaliascnt{remark}{theorem}
\newtheorem{remark}[remark]{Remark}
\aliascntresetthe{remark}

\newaliascnt{hypothesis}{theorem}

\aliascntresetthe{hypothesis}

\newaliascnt{property}{theorem}

\aliascntresetthe{property}

\usepackage[utf8]{inputenc}
\usepackage{amsmath}
\usepackage{enumerate}
\usepackage{bm}
\usepackage[T1]{fontenc}
\usepackage{amsthm}
\usepackage{amsfonts}
\usepackage{amssymb}
\usepackage{graphicx}
\usepackage{enumerate}
\usepackage[english]{babel}
\usepackage{comment}
\usepackage{amssymb}
\usepackage{amsthm}
\usepackage{xcolor}
\usepackage{dutchcal}
\usepackage{esint}

\usepackage[title]{appendix}
\usepackage{xcolor}

\usepackage{amssymb,amsmath,amsthm}
\usepackage{hyperref}
\usepackage{xcolor}
\usepackage{graphicx}
\usepackage[title]{appendix}
\usepackage{scalerel}[2014/03/10]
\usepackage{stackengine}

\def\dis{d}

\def\fp{{-\Delta_{p}}}

\def\fps{(-\Delta)_{p}^{s}}
\def\fqs{(-\Delta)_{q}^{s}}

\def\ueps{{u_\varepsilon}}
\def\veps{{v_\varepsilon}}

\def\pst{{p_{*}}}
\def\qst{{q_{s}^{*}}}
\def\jl{{J_{\lambda}}}
\def\Jl{{\mathcal{J}_{\lambda}}}
\def\Nl{\mathcal{N}_{\lambda}}
\def\nl{N_{\lambda}}
\def\mpqcr{{m(N,p,q,s)}}
\def\anorm{{\|a\|_{L^{\frac{r}{r-\delta}}}}}
\def\bnorm{{\|b\|_{L^{\infty}}}}

\def\gammau{{\gamma_{u}}}
\def\pt{{p_{*}}}

\def\fps{(-\Delta)_{p}^{s}}
\def\fqs{(-\Delta)_{q}^{s}}

\newcommand{\RR}{\mathbb{R}}

\newcommand{\wop}{W^{1,p}_0}

\newcommand{\wsq}{W^{s,q}_0}

\newcommand{\wopo}{W^{1,p}_0(\Omega)}
\newcommand{\wsqo}{W^{s,q}_0(\Omega)}

\newcommand{\w}{\mathcal{W}}

\newcommand{\wo}{\mathcal{W}_0}

\newcommand{\vep}{v_{\varepsilon}}

\DeclareMathOperator*{\prpr}{\prime\prime}

\newenvironment{sketch}{%
  \proof}{\endproof}
\newcommand{\Addresses}{{
		\footnote{
				\footnotesize
\noindent \textsuperscript{1}Indian Institute of Science Education and Research, Thiruvananthapuram 695551, India  \par\nopagebreak 
   \noindent \textsuperscript{2}Universite de Pau et des Pays de l’Adour, LMAP (UMR E2S-UPPA CNRS 5142) Bat. IPRA, Avenue de l’Universite F-64013 Pau, France  \par\nopagebreak 
\noindent \textsuperscript{A}\textit{e-mail:} \texttt{dhanya.tr@iisertvm.ac.in}
\noindent \textsuperscript{B}\textit{e-mail:} \texttt{jacques.giacomoni@univ-pau.fr}
\noindent \textsuperscript{C}\textit{e-mail:} \texttt{ritabrata20@iisertvm.ac.in}

			 \noindent \textsuperscript{*}Corresponding author.

			\medskip\noindent
			{\bf Acknowledgments:} 
R. Dhanya was supported by SERB MATRICS grant MTR/2022/000780 when this work was being carried out. J. Giacomoni was partially funded by IFCAM (Indo-French Centre for Applied Mathematics), IRL CNRS 3494. Ritabrata Jana expresses gratitude for the financial assistance provided by the Prime Minister's Research Fellowship during the execution of this research. The authors express their sincere gratitude to Dr. Vladimir Bobkov for the valuable discussions on the generalized eigenvalue problem, which greatly contributed to the completion of the proof of the nonexistence result.		
}}}
\begin{document}
\title[Multiplicity Results for Mixed Local Nonlocal Equations  ]{ Multiplicity Results for Mixed Local Nonlocal Equations With Indefinite Concave-Convex Type Nonlinearity	\Addresses	}
	\author[ R. Dhanya ]
	{ R. Dhanya\textsuperscript{1,A}} 
    \author[Jacques Giacomoni]
	{Jacques Giacomoni\textsuperscript{2,B,*}}
 	\author[Ritabrata Jana]
	{Ritabrata Jana\textsuperscript{1,C}}

\maketitle
\begin{abstract}
In this article we examine the multiplicity of non-negative solutions to mixed local-nonlocal equations involving \(\fp + \fqs\) in a bounded smooth domain. The nonlinearity incorporates a parameter \(\lambda > 0\), a sublinear term, and a superlinear term, with sign-changing weight functions \(a(x)\) and \(b(x)\). Under suitable conditions, we establish the existence of at least two distinct nontrivial non-negative solutions in both the subcritical and critical regimes via fibering map analysis and constrained minimization on the Nehari manifold. Additionally, for \(p \not = q\), we obtain a nonexistence result for large \(\lambda\) by analyzing the associated generalized eigenvalue problem.
\end{abstract}
\\
\section{Introduction}
In this work, we consider the problem
\begin{equation}\label{generaleqn}
    \begin{aligned}
        \fp u + \fqs u &= \lambda \left(a(x)|u|^{\delta-2}u+b(x)|u|^{r-2}u\right) \text{ in } \Omega
        \text{ and }
        u=0 \text{ in } \Omega^{c},
    \end{aligned}
\end{equation}
where $sq<p<N,$  and $\Omega$ is a bounded $C^{1,1}$ domain in $\RR^{N}.$ The parameter $\lambda$ is a positive real  number. We assume the sublinearity condition on $\delta,$ given by  $1<\delta<\min\{p,q\}.$ For $ 1<p<N, $ the critical Sobolev exponent $\pst$ and the fractional critical exponent $\qst$ are defined as $\pst=\frac{Np}{N-p}$ and $\qst=\frac{Nq}{N-sq}.$  Furthermore, we assume that $\max\{p,q\}\leq r\leq \max \{\pst,\qst\},$ with  $a(x)\in L^{\frac{r}{r-\delta}}(\Omega)$ and $b(x)\in L^{\infty}(\Omega)$. The local part of the operator is defined as $ \fp u:= -\nabla.(|\nabla u|^{p-2}\nabla u).$ The nonlocal part, the fractional $q$-Laplacian, denoted by $\fqs$, is defined as 
\begin{equation*}
(-\Delta)_{q}^su(x)= 2\thinspace PV. \int_{\mathbb{R}^N}\frac{|u(x)-u(y)|^{q-2}(u(x)-u(y))}{|x-y|^{N+sq}} \,dy \qquad x\in \Omega,
\end{equation*}
up to a suitable normalization constant. 
These operators emerge from the interaction of two stochastic processes operating at different scales: a classical random walk and a Lévy flight. When a particle transitions between these processes according to a defined probability distribution, the resulting limit diffusion equation is governed by a mixed local-nonlocal  operator. A detailed discussion of this phenomenon can be found in the appendix of \cite{DPV23}. These operators are also significant in various applications, including biological sciences, as highlighted in \cite{DV21} and related references, and in the study of heat transport in magnetized plasmas, as explored in \cite{BD13}.
   The qualitative properties such as regularity and comparison principles of solutions to the equation \(\fp u + \fqs u = f(x,u)\) subject to zero Dirichlet boundary conditions have been extensively analyzed in the homogeneous setting (\(p = q\)), as documented in \cite{BDVV22, CKSR12, GK22, GL23a, SVWZ22} among others. 
   More recently, for the non-homogeneous case 
   $(p\not = q)$, studies \cite{AC23, FM22, DGJ24} have focused on establishing interior and boundary $C^{1,\alpha}$ regularity results.
    \par
The existence theory for the solutions of Laplace equation with critical nonlinearity with a perturbation dates back to the seminal work of 
\cite{BN83}. Their work introduced a novel variational approach to address the difficulties associated with the critical Sobolev exponent \(2^*\), ultimately proving the existence of nonnegative solutions. Subsequently, \cite{ABC94} advanced this study by establishing existence and multiplicity results for the concave-convex type nonlinearity. Additionally, \cite{ABC94} provided nonexistence results for sufficiently large parameters, proved the existence of a minimal solution for small parameters, and analyzed the asymptotic behavior of the \(L^\infty\) norm of the minimal solution as the parameter approached zero. A similar problem for \(p\)-Laplacian was later explored in \cite{AGP96}. For further developments on multiplicity results in problems involving indefinite nonlinearities and local operators, we refer to \cite{FGU09, SCGG19} and the related references.  
   In the context of nonlocal linear operators, \cite{BCS13} investigated a concave-convex type nonlinearity. Their work provided a complete characterization of the parameter ranges ensuring the existence of solutions and established a multiplicity result. The influence of the critical exponent for the fractional Laplacian \( (-\Delta)^s \) was later addressed in \cite{CPS14}. Further developments on nonlocal elliptic equations with concave-convex nonlinearities were presented in \cite{SW15}, where the existence of at least six distinct solutions was demonstrated. Chen et al. \cite{CD15} employed fibering maps and the Nehari manifold approach to obtain multiple solutions for a related problem involving a nonlocal operator. For results concerning nonlocal and nonlinear operators with indefinite nonlinearities, we cite \cite{GKS20, BM19} and the associated references.
   \par
Returning to mixed local-nonlocal operators, a Brezis-Nirenberg-type result for the mixed local-nonlocal linear operator (i.e., \( p = q = 2 \)) was explored in \cite{BDVV22b}. The necessary and sufficient conditions for the existence of a unique positive weak solution to certain sublinear Dirichlet problems governed by a mixed local-nonlocal quasilinear operator (i.e., when \( p = q \)) were established in \cite{BMV24}. In a related direction, \cite{DFV24} employed a combination of variational and topological methods to study the existence and multiplicity of nontrivial solutions to problems driven by \( \fp + \fps \), where the nonlinearity is given by \( f(x,u) = \lambda |u|^{q-2} u + |u|^{\pt-2} u \), with \( \lambda \) being a real parameter and \( q \) possibly sublinear, linear, or subcritical. Furthermore, \cite{BV24} recently proved the existence of at least two positive weak solutions for a singular and critical semilinear elliptic problem involving a mixed local-nonlocal operator, in the spirit of \cite{Hai03}. A mixed local-nonlocal critical semilinear elliptic problem with a sublinear perturbation, given by \( -\Delta u + \varepsilon(-\Delta)^s u = \lambda u^p + u^{2_* - 1} \) under zero Dirichlet boundary conditions, where \( 0 < p < 1 \), \( \varepsilon \in (0,1] \), and \( \lambda > 0 \), was analyzed in \cite{BV24b}. Their work establishes the existence of a second positive weak solution for a concave-convex type semilinear mixed local-nonlocal problem, following the spirit of \cite{ABC94}.  
\par
As discussed in \cite{FGU09, Pai11} for local operators and \cite{GKS20} for nonlocal ones, the presence of sign-changing weight functions combined with subcritical and critical nonlinearities introduces significant challenges. Even in the purely local case, this complicates the application of established tools such as the nonquadraticity condition by Costa–Magalhães \cite{CM94}. Notably, in the context of mixed local-nonlocal operators, including the case \(p = q \), sign-changing weight functions have not yet been explored. To bridge this gap, we examine the multiplicity of solutions for quasilinear elliptic equations involving the nonhomogeneous mixed local-nonlocal operator \(\fp + \fqs\) with indefinite concave-convex nonlinearities, where the weight functions \(a(x)\) and \(b(x)\) are permitted to change sign. The nonhomogeneous nature of the operator further complicates the analysis, rendering the investigation of multiplicity and nonexistence results novel, even for constant weight functions \(a(x) \equiv b(x) \equiv 1\). The existing literature on mixed local-nonlocal operators with critical nonlinearity, such as \cite{BV24, DFV24}, focuses only on homogeneous operators. This homogeneity enables the use of Talenti functions associated with the p-Laplacian to establish that the associated  functional verifies the Palais Smale condition when energy level is beyond a certain threshold.
In contrast, our case requires estimating the \(\wsq\) norm of the Talenti functions of the \(\fp\), which, to the best of our knowledge, has not been previously addressed.
\par
Using fibering map analysis and constrained minimization on specific subsets of the Nehari manifold, we establish the existence of at least two distinct nontrivial non-negative solutions to \eqref{generaleqn} for sufficiently small values of \(\lambda\). We first analyze the case \(p < q\), where both operators significantly influence the behavior of the solutions. Here, we consider a subcritical nonlinearity with a sublinear perturbation, allowing \(\lambda\) to appear in both terms while permitting the coefficients \(a\) and \(b\) to change sign. Furthermore, we prove a nonexistence result for sufficiently large \(\lambda\) when \(p < q\), even when \(a\) and \(b\) are positive. To prove this result, we derive a few key properties of the threshold curve in the parameter plane of the generalized eigenvalue problem. The novel observations we provide related to the generalized eigenvalue problem  of mixed local-nonlocal operator may be of interest to researchers from different perspectives.  
\par
Finally, we examine the critical perturbation case when \(q \leq p\), where the local operator dominates. For sufficiently small \(\lambda\), we show the existence of two nonnegative solutions in the presence of both nonlinearities, assuming \(a(x)\) is continuous and \(b(x) \equiv 1\). In this critical regime, the Palais-Smale condition fails globally for the associated energy functional. However, we prove that it holds for energy levels below the first critical level. By employing Talenti functions for the local operator, we demonstrate that the energy remains below this critical threshold, leading to the multiplicity result.
Finally, we explore a Brezis-Nirenberg type problem involving critical perturbations when \(q \leq p\), with \(\lambda\) appearing only in the sublinear term. Using the asymptotic behavior of Talenti functions, we establish a multiplicity result for a specific range of sublinearity. The main results of this work are summarized below.
   \begin{theorem}\label{subcrmul}
    Let \( p \leq q \) and \( r < \max\{\pt,q_{s}^{*}\} \). Additionally, suppose that $a(x)\in L^{\frac{r}{r-\delta}}(\Omega)$ and $b(x)\in L^{\infty}(\Omega)$. Then, there exists \( \lambda_0 > 0 \) such that for all \( \lambda \in (0, \lambda_0) \), the problem \eqref{generaleqn} admits at least two nonnegative, non-trivial solutions.
\end{theorem}
\begin{theorem}\label{nonexsubcr}
    Suppose $\inf a = \alpha > 0,$ $a(x)\in L^{\frac{r}{r-\delta}}(\Omega),$ $\inf b = \beta > 0,$   and $b(x)\in L^{\infty}(\Omega)$. Assume that $sq<p,$ $p\neq q$ and $r < \max\{\pt,q_{s}^{*}\}$, then there exists a threshold $\Lambda_* > 0$ such that for all $\lambda > \Lambda_*$, problem \eqref{generaleqn} admits only the trivial solution.
\end{theorem}
\begin{theorem}\label{multsolcri}
    Let \( q \leq p \) and $r = p_*$ in  \eqref{generaleqn}. Additionally, assume that $b(x) \equiv 1$ and $a(x)\in L^{\infty}(\Omega)$ is a continuous function. There exists a constant $\Lambda_{0} > 0$ such that for all $\lambda \in (0, \Lambda_{0})$, equation \eqref{generaleqn} admits at least two distinct, nontrivial, nonnegative solutions for any $\delta < q$.
\end{theorem}
\begin{theorem}\label{bnthm}
   Suppose that \( q \leq p \) and \( r = \pt \) in equation \eqref{generaleqn} and define \[ \mpqcr = \min \left\{ \frac{q(N-p)}{p(p-1)}, q(1-s) + N \left( 1 - \frac{q}{p} \right) \right\}. \] 
   Assume that one of the following conditions holds:  
    \begin{equation*}
        \max\left\{\frac{Np}{\mpqcr+N-p},\pt\left(1-\frac{1}{p}\right)\right\} < \delta < q,
    \end{equation*}  
    or  
    \begin{equation*}
        \delta < \min\left\{ q, \pt\left(1-\frac{1}{p}\right)\right\} \quad \text{and} \quad 0<s < 1 - \frac{1}{q} \left(\frac{N-p}{p-1} - N\left(1 - \frac{q}{p}\right) \right).
    \end{equation*}  
    Additionally, assume that \( b(x) \equiv \lambda^{-1} \) and that \( a(x)\in L^{\infty}(\Omega) \) is a continuous function. Then, there exists a constant \( \bar\Lambda_{00} > 0 \) such that for all \( \lambda \in (0, \bar\Lambda_{00}) \), equation \eqref{generaleqn} has at least two distinct, nontrivial, nonnegative solutions.  
\end{theorem}
The article is organized as follows: Section \ref{prlm} introduces the necessary definitions and notations that form the foundation for the subsequent analysis. Section \ref{nehari} presents a detailed study of the Nehari manifold and the fibering map analysis. In Section \ref{pscnd}, we get the details about the Palais–Smale condition .  In Section \ref{subcritical}, we establish the existence of two solutions for the cases \( p \leq q \) and \( r < \max\{\pt, \qst\} \). Section \ref{nonexistence} examines the nonexistence results under the conditions \( p \leq q \), \( r < \max\{\pt, \qst\} \), and when the weight functions are positive.  Section \ref{critical} establishes a multiplicity result for \( q \leq p \) and \( r = \pt \). Finally, in Section \ref{brezisnirenberg}, we investigate a Brezis-Nirenberg-type problem in the setting where \( b(x) \equiv \lambda^{-1} \), \( q \leq p \), and \( r = \pt \).
\par
\textbf{Notations:}
Throughout this article, unless stated otherwise, the symbols \( k, M, C \), etc., represent generic positive constants, whose values may vary even within the same line. We assume that \( p, q > 1 \) and \( p > sq \) throughout our analysis. Given any \( a \in \mathbb{R} \), we define \( a_{+} := \max\{a,0\} \). Additionally, for any \( a \in \mathbb{R} \) and \( t > 0 \), we use the notation \([a]^t := |a|^{t-1}a\). The open ball of radius \( R > 0 \) centered at \( x_0 \in \mathbb{R}^N \) is denoted by \( B_R(x_0) \), and when the center is not relevant, we omit its notation. For \( t > 1 \), we define the \( L^t \)-norm as \( \|\cdot\|_{t} \equiv \|\cdot\|_{L^t} \). Given a subset \( S \subset \mathbb{R}^{2N} \), we introduce the following functionals:  
\begin{equation*}
    \begin{aligned}
    \mathfrak{A}_{t}(u,v) &= \int_{\Omega} |\nabla u|^{t-2}\nabla u \cdot \nabla v \,dx, \\
    A_{t}(u,v,S) &= \int_{S} \frac{|u(x)-u(y)|^{t-2}(u(x)-u(y))(v(x)-v(y))}{|x-y|^{N+ts}} \, dx \, dy.
    \end{aligned}
\end{equation*}  
We use the notation \( f(x) = O(g(x)) \) to indicate that there exists a positive constant \( M \) and a real number \( x_0 \) such that  
\begin{equation*}
    |f(x)| \leq M |g(x)| \quad \text{for all } x \geq x_0.
\end{equation*}  
Similarly, we write \( f(x) = o(g(x)) \) if  
\begin{equation*}
    \lim_{x\rightarrow\infty} \frac{|f(x)|}{|g(x)|} = 0.
\end{equation*}
We say $f\lesssim g$  in $\Omega$ if there exists a $C(N,s,p,q,\Omega)>0$ such that $f\leq C g$ in $\Omega.$ 
 
\section{Preliminaries}\label{prlm}
We recall that for $E \subset \mathbb{R}^N$, the Lebesgue space $L^t(E)$, $1 \leq t < \infty$, is defined as the space of $t$-integrable functions $u: E \to \mathbb{R}$ with the finite norm 
\begin{equation*}
    \|u\|_{L^t(E)} = \left( \int_E |u(x)|^t \, dx \right)^{1/t}.
\end{equation*}
 The Sobolev space $W^{1,t}(\Omega)$, for $1 \leq t < \infty$, is defined as the Banach space of locally integrable weakly differentiable functions $u: \Omega \to \mathbb{R}$ equipped with the following norm 
\begin{equation*}
    \|u\|_{W^{1,t}(\Omega)} = \|u\|_{L^t(\Omega)} + \|\nabla u\|_{L^t(\Omega)}.
\end{equation*}
The space $\mathfrak{W}^{1,t}_0(\Omega)$ is defined as the closure of the space $C_c^\infty(\Omega)$ in the norm of the Sobolev space $W^{1,t}(\Omega)$, where $C_c^\infty(\Omega)$ is the set of all smooth functions whose supports are compactly contained in $\Omega$. The space \( W^{1,t}_0(\Omega) \) is defined as the set of all functions in \( \mathfrak{W}^{1,t}_0(\Omega) \) that vanish outside \( \Omega \). For a measurable function $u:\mathbb{R}^{N}\rightarrow \mathbb{R}$, we define Gagliardo seminorm
\begin{equation*}
	[u]_{s,t}:=[u]_{W^{s,t}(\mathbb{R}^N)}:=  \left(\int_{\mathbb{R}^N \times \mathbb{R}^N} \dfrac{|u(x)-u(y)|^{t}}{|x-y|^{N+st}}\,dx\, dy\right)^{1/t},
\end{equation*}
for $1<t<\infty$ and $0<s<1.$
We consider the space $W^{s,t}(\mathbb{R}^N)$ defined as  
\begin{equation*}
	W^{s,t}(\mathbb{R}^{N}):= \left \{u \in {L}^{t}(\mathbb{R}^{N}):[u]_{s,t}<\infty  \right\}.
\end{equation*}
The space $W^{s,t}(\mathbb{R}^{N})$ is a Banach space with respect to the norm 
\begin{equation*}
\|u\|_{{W^{s,t}(\mathbb{R}^{N})}}=\left( \|u\|^{t}_{L^{t}(\mathbb{R}^{N})} + [u]^{t}_{{W^{s,t}(\mathbb{R}^{N})}}\right)^{\frac{1}{t}} . 	
\end{equation*}
A comprehensive examination of the fractional Sobolev Space and its properties are presented in \cite{DPV12}. To address the Dirichlet boundary condition, we naturally consider the space $W^{s,t}_0(\Omega)$ defined as
\begin{equation*}
	W_{0}^{s,t}(\Omega):= \left \{u \in W^{s,t}(\mathbb{R}^{N}):u=0\medspace\text{in}\medspace \mathbb{R}^{N}\setminus \Omega \right\}.
\end{equation*}
This is a separable, uniformly convex Banach space endowed with the norm $\|u\|=	\|u\|_{{W^{s,t}(\mathbb{R}^{N})}}$.
For \( N > t \), the critical Sobolev exponent is given by \( t^* = \frac{Nt}{N - t} \). A fundamental embedding result states that for any bounded open subset \( \Omega \) of class \( C^1 \) in \( \mathbb{R}^N \), there exists a constant \( C \equiv C(N, \Omega) > 0 \) such that for all \( u \in C_c^{\infty}(\Omega) \), the inequality  
\begin{equation*}
    \|u\|_{L^{t^*}(\Omega)} \leq C \int_{\Omega} |\nabla u|^t \, dx
\end{equation*}  
holds. The best constant in the Sobolev embedding for \( W_0^{1,p}(\Omega) \) is defined as  
\begin{equation*}
    S_{p}:=\underset{\underset{v\not\equiv 0}{v\in W_0^{1,p}(\Omega)}}{\inf} \frac{\|v\|_{W^{1,p}(\Omega)}^{p}}{\|v\|_{L^\pt(\Omega)}^{p}}.
\end{equation*}   
In general in this work we define 
\begin{equation*}
    \begin{aligned}
        S_{rp}:=\underset{\underset{v\not\equiv 0}{v\in W_0^{1,p}(\Omega)}}{\inf} \frac{\|v\|_{W^{1,p}(\Omega)}^{p}}{\|v\|_{L^r(\Omega)}^{p}}. 
    \end{aligned}
\end{equation*} 
Moreover, the inclusion map  
\begin{equation*}
    W_0^{1,t}(\Omega) \hookrightarrow L^r(\Omega)
\end{equation*}  
is continuous for \( 1 \leq r \leq t^* \), and compact except when \( r = t^* \). Similarly, the embedding  
\begin{equation*}
    W^{s,t}_{0}(\Omega)\hookrightarrow L^{r}(\Omega)
\end{equation*}  
is continuous for \( 1\leq r\leq t^{*}_{s} := \frac{Nt}{N-ts} \) and compact for \( 1\leq r < t^{*}_{s} \). Due to the continuous embedding of \( W^{s,t}_{0}(\Omega) \) into \( L^{r}(\Omega) \) for \( 1\leq r\leq t^{*}_{s} \), we define an equivalent norm on \( W^{s,t}_{0}(\Omega) \) as  
\begin{equation*}
	\|u\|_{W^{s,t}_{0}}:=\left(\int_{\mathbb{R}^N \times \mathbb{R}^N} \dfrac{|u(x)-u(y)|^{t}}{|x-y|^{N+st}}\,dx\, dy\right)^{1/t}.
\end{equation*}  
The best constants for the embeddings  \( W_0^{s,q}(\Omega) \) are similarly given by  
\begin{equation*}
    S_{q}:=\underset{\underset{v\not\equiv 0}{v\in W_0^{s,q}(\Omega)}}{\inf} \frac{\|v\|_{W^{s,q}(\Omega)}^{q}}{\|v\|_{L^\qst(\Omega)}^{q}}\quad\text{ and }\quad S_{rq}:=\underset{\underset{v\not\equiv 0}{v\in W_0^{s,q}(\Omega)}}{\inf} \frac{\|v\|_{W^{s,q}(\Omega)}^{q}}{\|v\|_{L^r(\Omega)}^{q}} .
\end{equation*}  
The dual space of \( W^{s,t}_{0}(\Omega) \) is denoted by \( W^{-s,t'}(\Omega) \) for \( 1<t<\infty \). 
There may not be a continuous embedding between \( W^{1,p}(\Omega) \) and \( W^{s,q}(\Omega) \) if $p<q$. To address this, when analyzing weak solutions associated with the operator \( \fp + \fqs \), for $p<q$ we consider the space  
\begin{equation*}
    \mathcal{W}(\Omega) = W^{1,p}(\Omega) \cap W^{s,q}(\Omega),
\end{equation*}  
equipped with the norm \( \|\cdot\|_{\mathcal{W}(\Omega)} = \|\cdot\|_{W^{1,p}(\Omega)} + \|\cdot\|_{W^{s,q}(\Omega)} \). To incorporate the zero Dirichlet boundary condition, we define  
\begin{equation*}
    \mathcal{W}_0(\Omega) = W_0^{1,p}(\Omega) \cap W_0^{s,q}(\Omega).
\end{equation*}  
Finally, the dual space of \( \mathcal{W}_0(\Omega) \) is denoted by \( \mathcal{W}^{\prime}(\Omega) \).
\begin{definition}
    We say that \( u \in \mathcal{W}_0(\Omega) \) is a subsolution (supersolution) to the problem \eqref{generaleqn} if for every non-negative test function \( \varphi \in \mathcal{W}_0(\Omega) \), the following inequality holds:
    \begin{equation*}
        \begin{aligned}
            \int_{\Omega} |\nabla u|^{p-2} \nabla u \cdot \nabla \varphi &+ \int_{\RR^N} \int_{\RR^N} \frac{|u(x) - u(y)|^{q-2} (u(x) - u(y)) (\varphi(x) - \varphi(y))}{|x - y|^{N+sq}} \,dx \,dy
           \\
           &\leq (\geq) \ \int_{\Omega} \lambda \left( a(x) |u|^{\delta-1} u + b(x) |u|^{r-1} u \right) \varphi \,dx.
        \end{aligned}
    \end{equation*}
    A function \( u \) is said to be a solution if it satisfies both the subsolution and supersolution conditions.
\end{definition}
The associated energy functional \( J_{\lambda} : \mathcal{W}_0(\Omega) \to \mathbb{R} \) is given by
\begin{equation*}
    J_{\lambda}(u) := \frac{1}{p} \|u\|_{\wopo}^{p} + \frac{1}{q} \|u\|_{\wsqo}^{q} - \lambda \int_{\Omega} \left( \frac{a(x)}{\delta} |u|^{\delta} + \frac{b(x)}{r} |u|^{r} \right) dx.
\end{equation*}
We define \( \lambda \) as an eigenvalue of \( \fp \) if there exists a nontrivial solution to  \( \fp u = \lambda |u|^{p-2} u \) with zero Dirichlet boundary conditions. Similarly, \( \lambda \) is an eigenvalue of \( \fqs \) if there exists a nontrivial solution to the problem \( \fqs u = \lambda |u|^{q-2} u \) with zero Dirichlet boundary conditions. The first positive eigenvalues \( \lambda_{1p} \) and \( \lambda_{1q} \) are obtained by minimizing the Rayleigh quotient:
\begin{equation*}
    \lambda_{1p} := \underset{\underset{\|u\|_{L^p} > 0}{v \in \wopo}}{\inf} \frac{\|v\|_{\wop}^{p}}{\|v\|_{L^p}^{p}}, \quad
    \lambda_{1q} := \underset{\underset{\|u\|_{L^q} > 0}{v \in \wsqo}}{\inf} \frac{\|v\|_{\wsq}^{q}}{\|v\|_{L^q}^{q}}.
\end{equation*}
By \cite[Proposition 2.1]{BF22} and \cite{BT15}, the eigenfunction \( \phi_q \) corresponding to \( \lambda_{1q} \) and the eigenfunction \( \phi_p \) corresponding to \( \lambda_{1p} \) have a constant sign. Moreover, the eigenvalues \( \lambda_{1q} \) and \( \lambda_{1p} \) are simple and isolated.
\\
Consider the linear space  
\[
C^1_0(\bar\Omega) := \{ u \in C^1(\bar\Omega) : u|_{\Omega^c} = 0 \},
\]
which is a Banach space under the standard \( C^1 \)-norm. We define its positive cone as  
\[
C^+ := \{ u \in C^1_0(\Omega) : u(x) \geq 0 \text{ for all } x \in \Omega \}.
\]
The interior of \( C^+ \) is nonempty and given by  
\[
\operatorname{int}(C^+) = \left\{ u \in C^+ : u(x) > 0 \text{ for all } x \in \Omega, \quad \frac{\partial u}{\partial n}(x) < 0 \text{ for all } x \in \partial\Omega \right\}.
\]
\section{Nehari Manifold and Fibering Map Analysis}\label{nehari}
\subsection{The Nehari Manifold}
The primary objective of this section is to analyze the critical points of the fibering maps associated with the energy functional \( J_{\lambda} \). For a comprehensive discussion on the Nehari method, we refer the reader to \cite{BZ03,Wil96}. The Nehari manifold corresponding to \( J_{\lambda} \) is defined as  
\begin{equation*}
    \begin{aligned}
        N_{\lambda} &:= \{u \in \mathcal{W}_0(\Omega) \setminus \{0\} : \langle J_{\lambda}^{\prime}(u), u \rangle = 0\} \\
        &= \{u \in \mathcal{W}_0(\Omega) \setminus \{0\} : \|u\|_{\mathcal{W}_0} = \lambda \int_{\Omega} (a(x)|u|^{\delta} + b(x)|u|^{r}) \,dx\}.
    \end{aligned}
\end{equation*}
where \( \langle \cdot, \cdot \rangle \) denotes the duality pairing between \( \mathcal{W}_0 \) and its dual space \( \mathcal{W}^{\prime} \). Since the mapping \( u \mapsto \langle J_{\lambda}^{\prime}(u), u \rangle \) is a \( C^{1} \) functional, it follows that \( N_{\lambda} \) forms a \( C^{1} \) submanifold of \( \mathcal{W}_0(\Omega) \). Furthermore, every solution of problem \eqref{generaleqn} belongs to \( N_{\lambda} \). As an initial step, we demonstrate that \( J_{\lambda} \) is coercive and bounded from below on \( N_{\lambda} \), which enables us to obtain a ground state solution for problem \eqref{generaleqn}.

\begin{lemma}
The functional $\jl$ is coercive and bounded from below on $\nl.$
\end{lemma}
\begin{proof}
    Since $u\in \nl,$ $\|u\|_{\wopo}+ \|u\|_{\wsqo}-\lambda \int_{\Omega} a(x)|u|^{\delta}=\lambda \int_{\Omega} b(x)|u|^{r}.$ Thus using the Holder inequality we have 
    \begin{equation*}
\begin{aligned}
    \jl(u)\geq \left(\frac{1}{p}-\frac{1}{r}\right) \|u\|_{\wop}^{p} + \left(\frac{1}{q}-\frac{1}{r}\right)\|u\|_{\wsq}^{q}-\lambda\left(\frac{1}{\delta}-\frac{1}{r}\right) \|a\|_{L^{\frac{r}{r-\delta}}} S_{p}^{\frac{-\delta}{p}}\|u\|_{\wop}^{\delta}.
\end{aligned}
    \end{equation*}
   Since $\delta<\min\{p,q\},$ we have $\jl$ is coercive and bounded from below on $\nl.$
\end{proof}
We define the fibering map associated with \( J_{\lambda} \) as \( \gamma_{u}:\mathbb{R}_{+}\to\mathbb{R} \) by  
\begin{equation}\label{fibermap}
    \gamma_u(t) = J_{\lambda}(tu).
\end{equation}
Explicitly, we set  
\begin{equation*}
    \begin{aligned}
        \gamma_u(t) = \frac{t^p}{p} \|u\|_{\wop}^p + \frac{t^q}{q} \|u\|_{\wsq}^q - \lambda \int_{\Omega} \left(\frac{t^{\delta}}{\delta} a(x) |u|^{\delta} + \frac{t^r}{r} b(x) |u|^r \right)dx.
    \end{aligned}
\end{equation*}
Fibering maps are widely studied alongside the Nehari manifold to establish the existence of critical points for \( J_{\lambda} \). In particular, for problems involving concave–convex nonlinearities, it is crucial to analyze the geometry of \( \gamma_u \) (see \cite{BZ03}). We note that \( \gamma_u \) is a \( C^1 \) function, and its derivative is given by  
\begin{equation}\label{firstderfibermap}
    \begin{aligned}
        \gamma_u^{\prime}(t) = t^{p-1} \|u\|_{\wop}^p + t^{q-1} \|u\|_{\wsq}^q - \lambda \int_{\Omega} \left( t^{\delta-1} a(x) |u|^{\delta} + t^{r-1} b(x) |u|^r \right)dx.
    \end{aligned}
\end{equation}
It follows that \( tu \in N_{\lambda} \) if and only if \( \gamma_u^{\prime}(t) = 0 \), and in particular, \( u \in N_{\lambda} \) if and only if \( \gamma_u^{\prime}(1) = 0 \). This observation implies that locating stationary points of the fibering map \( \gamma_u \) is sufficient to identify critical points of \( J_{\lambda} \) on \( N_{\lambda} \). Moreover, \( \gamma_u \) is twice differentiable function, and its second derivative is given by  
\begin{equation}\label{secondderfibermap}
    \begin{aligned}
        \gamma_u^{\prime\prime}(t) = (p-1)t^{p-2} \|u\|_{\wop}^p &+ (q-1)t^{q-2} \|u\|_{\wsq}^q \\
        &- \lambda \int_{\Omega} \left( (\delta-1) t^{\delta-2} a(x) |u|^{\delta} + (r-1) t^{r-2} b(x) |u|^r \right)dx.
    \end{aligned}
\end{equation}
Naturally, we classify \( N_{\lambda} \) into three subsets corresponding to local minima, local maxima, and points of inflection:  
\begin{equation*}
    \begin{aligned}
        N_{\lambda}^{+} &= \{u \in N_{\lambda} : \gamma_u^{\prime\prime}(1) > 0\}, \\
        N_{\lambda}^{-} &= \{u \in N_{\lambda} : \gamma_u^{\prime\prime}(1) < 0\}, \\
        N_{\lambda}^{0} &= \{u \in N_{\lambda} : \gamma_u^{\prime\prime}(1) = 0\}.
    \end{aligned}
\end{equation*}
Furthermore, we define the following critical levels:
\begin{equation*}
    \theta_{\lambda} := \inf_{u \in N_{\lambda}} J_{\lambda}(u), \quad 
    \theta_{\lambda}^{\pm} := \inf_{u \in N_{\lambda}^{\pm}} J_{\lambda}(u).
\end{equation*}
Following the classical work of Drábek and Pohozaev \cite{DP97}, we conclude that  if $\nl^0=\varnothing$ then any minimizer of \( J_{\lambda} \) on the Nehari manifold \( N_{\lambda} \)  is a critical point of \( J_{\lambda} \) in the entire space \( \mathcal{W}_0(\Omega) \). In particular, we establish the following lemma.
\begin{lemma}
    If \( u \) is a minimizer of \( J_{\lambda} \) on \( N_{\lambda} \) and \( u \notin N_{\lambda}^{0} \), then \( u \) is a critical point of \( J_{\lambda} \).
\end{lemma}

We will now establish that, for sufficiently small values of \( \lambda \), the set \( N_{\lambda}^{0} \) is empty.  

\begin{lemma}\label{nlambdazero}
    There exists \( \lambda_{0} > 0 \) such that \( N_{\lambda}^{0} = \varnothing \) for all \( \lambda \in (0, \lambda_{0}) \).
\end{lemma}
\begin{proof}
    We first consider the case $u\in\nl$ and $\int_{\Omega}a(x)|u|^{\delta}=0.$ Since $u\in \nl,$ $\|u\|_{\wo}=\lambda \int_{\Omega} b(x)|u|^{r},$ and 
    \begin{equation*}
        \begin{aligned}
            \gammau^{\prime\prime}(1)
            &=(p-r)\|u\|_{\wop}+ (q-r) \|u\|_{\wsq}<0.
        \end{aligned}
    \end{equation*}
    Next consider the case $u\in\nl$ and $\int_{\Omega}a(x)|u|^{\delta}\neq 0.$ Since $\gammau^{\prime}(1)=0,$  
    \begin{equation*}
        \begin{aligned}
            \lambda \int_{\Omega}a(x)|u|^{\delta} &= \|u\|_{\wop}^{p}+ \|u\|_{\wsq}^{q}-\lambda \int_{\Omega}  b(x)|u|^{r},
            \\
           \lambda \int_{\Omega}  b(x)|u|^{r}   &= \|u\|_{\wop}^{p}+ \|u\|_{\wsq}^{q}- \lambda \int_{\Omega}a(x)|u|^{\delta}.
        \end{aligned}
    \end{equation*}
    If $u\in\nl^{0}$ we have 
    \begin{equation}\label{abev}
        \begin{aligned}
       (p-\delta) \|u\|_{\wop}^{p}+ (q-\delta)\|u\|_{\wsq}^{q}-\lambda(r-\delta) \int_{\Omega}  b(x)|u|^{r}=0,
        \end{aligned}
    \end{equation}
     \begin{equation}\label{bbev}
        \begin{aligned}
       (r-p) \|u\|_{\wop}^{p}+ (r-q)\|u\|_{\wsq}^{q}-\lambda(r-\delta) \int_{\Omega}  a(x)|u|^{\delta}=0.
        \end{aligned}
    \end{equation}
    Define the functional \( E_{\lambda}: N_{\lambda} \to \mathbb{R} \) as  
\begin{equation*}
    E_{\lambda}(u) := \frac{r-p}{r-\delta} \|u\|_{\wop}^{p} + \frac{r-q}{r-\delta} \|u\|_{\wsq}^{q} - \lambda \int_{\Omega} a(x) |u|^{\delta}.
\end{equation*}
From \eqref{bbev}, it follows that \( E_{\lambda}(u) = 0 \) for all \( u \in N_{\lambda}^{0} \). Moreover, we obtain the lower bound  
\begin{equation*}
    E_{\lambda}(u) \geq \frac{r-p}{r-\delta} \|u\|_{\wop}^{p} - \lambda \anorm S_{rp}^{-\frac{\delta}{p}} \|u\|_{\wop}^{\delta}.
\end{equation*}
From \eqref{abev}, we have  
\begin{equation*}
    (p-\delta) \|u\|_{\wop}^{p} \leq \lambda (r-\delta) \int_{\Omega} b(x) |u|^{r}.
\end{equation*}
Applying Hölder’s inequality, we derive  
\begin{equation*}
    \left(\frac{p-\delta}{r-\delta} \frac{S_{rp}^{r/p}}{\lambda \|b\|_{L^{\infty}(\Omega)}} \right)^{\frac{1}{r-p}} \leq \|u\|_{\wop(\Omega)},
\end{equation*}
which leads to  
\begin{equation*}
    E_{\lambda}(u) \geq \|u\|_{\wop}^{\delta} \left( \frac{r-p}{r-\delta} \left(\frac{p-\delta}{r-\delta} \frac{S_{rp}^{r/p}}{\lambda \|b\|_{L^{\infty}(\Omega)}} \right)^{\frac{p-\delta}{r-p}} - \lambda \anorm S_{rp}^{-\frac{\delta}{p}} \right).
\end{equation*}
Thus, choosing  
\begin{equation*}
    \lambda_{0} = \underset{t=p,q}{\min}\left\{ \left(\frac{r-t}{r-\delta} \right)^{\frac{r-t}{r-\delta}} \left( \frac{(t-\delta) S_{rt}^{r/t}}{(r-\delta) \|b\|_{L^{\infty}(\Omega)}} \right)^{\frac{t-\delta}{r-\delta}} \left( \frac{S_{rt}^{\delta/t}}{\anorm} \right)^{\frac{r-t}{r-\delta}}\right\} > 0,
\end{equation*}
for sufficiently small \( \lambda \in (0, \lambda_{0}) \), we ensure that \( E_{\lambda}(u) > 0 \) for all \( u \in N_{\lambda}^{0} \), contradicting the fact that \( E_{\lambda}(u) = 0 \) in \( N_{\lambda}^{0} \). This completes the proof.
\end{proof}
\subsection{Fibering Map Analysis}\label{fibermapsec}
We now provide a complete characterization of the geometry of the fibering maps associated with problem \eqref{generaleqn}.
\begin{lemma}\label{fiberingmapanalysis}
\setlength\itemsep{1em}
    Let $u\in\wo\setminus\{0\}$ is a fixed function. Then ,
    \begin{enumerate}
        \item[(i)] Assume that $\int_{\Omega}b(x)|u|^{r}>0$ and $\int_{\Omega}a(x)|u|^{\delta}>0$ then there exists an unique $t_{\max}>0$ such that $\gamma^{\prpr}_{t_{\max}u}(1)=0.$ Moreover, there exist $ t_1(u,\lambda)< t_{\max} $ and $ t_2(u,\lambda)>t_{\max} $  such that \( t_{1}u \in \nl^{+} \) and \( t_{2}u \in \nl^{-} \) and  
        \[
\gammau^{\prime }(t) < 0 \text{ } \text{for all } t \in [0, t_{1}),\quad
\gammau^{\prime }(t) > 0 \text{ } \text{for all } t \in (t_{1}, t_{2}].
\]
\item[(ii)] Suppose \( \int_{\Omega} a(x) |u|^{\delta} < 0 \) and \( \int_{\Omega} b(x) |u|^{r} > 0 \) hold. Then there exists an unique $t_1(u,\lambda)>0$ such that $t_1 u\in\nl^{-}.$
\item[(iii)] There exists an unique $t_1(u,\lambda)>0$ such that $t_1 u\in\nl^{+}$ if \( \int_{\Omega} a(x) |u|^{\delta} > 0 \) and \( \int_{\Omega} b(x) |u|^{r} < 0 \).  
\item[(iv)] There is no critical point whenever  \( \int_{\Omega} a(x) |u|^{\delta} < 0 \) and \( \int_{\Omega} b(x) |u|^{r} < 0 \). 
    \end{enumerate}
\end{lemma}
\begin{proof}
We introduce the auxiliary $C^{1}$ function  $m_{u}:\RR_{+}\rightarrow\RR$ which is defined for a fixed $u\in\wo\setminus\{0\}$ as
\begin{equation*}
    \begin{aligned}
        m_{u}(t)= t^{(p-\delta)} \|u\|_{\wop}^{p}+ t^{(q-\delta)}\|u\|_{\wsq}^{q}-\lambda t^{(r-\delta)} \int_{\Omega}  b(x)|u|^{r}\text{ for } t\geq 0. 
    \end{aligned}
\end{equation*} 
Differentiating, we obtain
\begin{equation}
    \begin{aligned}
         m_{u}^{\prime}(t)=(p-\delta) t^{(p-\delta-1)} \|u\|_{\wop}^{p}+(q-\delta) t^{(q-\delta-1)}\|u\|_{\wsq}^{q}-\lambda (r-\delta)t^{(r-\delta-1)} \int_{\Omega}  b(x)|u|^{r}.
    \end{aligned}
\end{equation}
It follows that $tu\in\nl$ if and only if $t$ satisfies $m_{u}(t)=\lambda  \int_{\Omega}  a(x)|u|^{\delta}.$ Moreover if $tu\in\nl$ then the second derivative satisfies $\gamma_{tu}^{\prpr}(1)=t^{\delta+1}m_{u}^{\prime}(t).$ Now we analyse the behaviour of $\gammau$ based on the sign of $\int_{\Omega}a(x)|u|^{\delta}$ and $\int_{\Omega}b(x)|u|^{r}.$
\\
\textbf{Case $(i)$:}  $\int_{\Omega}a(x)|u|^{\delta}>0$ and $\int_{\Omega}b(x)|u|^{r}>0:$ We see $ m_{u}(t)\rightarrow -\infty$ as $t\rightarrow\infty,$ $ m_{u}(t)>0$ for $t$ small enough and $ m_{u}^{\prime}(t)<0$ for $t$ large enough. We claim that there exists unique $t_{max}>0$ such that $ m_{u}^{\prime}(t_{max})=0.$
\\
First we discuss the case $p>q.$ Then, we rewrite 
\begin{equation*}
    \begin{aligned}
          m_{u}^{\prime}(t)=t^{(q-\delta-1)}\left[(p-\delta) t^{(p-q)} \|u\|_{\wop}^{p}+(q-\delta) \|u\|_{\wsq}^{q}-\lambda (r-\delta)t^{(r-q)} \int_{\Omega}  b(x)|u|^{r}\right].
    \end{aligned}
\end{equation*}
Define the function  
\[
G_u(t) = (p-\delta) t^{(p-q)} \|u\|_{\wop}^{p} + (q-\delta) \|u\|_{\wsq}^{q} - \lambda (r-\delta)t^{(r-q)} \int_{\Omega} b(x)|u|^{r}.
\]
Our goal is to establish the existence of a unique \( t_{\max} > 0 \) such that \( G_u(t_{\max}) = 0 \). Consider  
\[
H_u(t) = -(p-\delta) t^{(p-q)} \|u\|_{\wop}^{p} + \lambda (r-\delta)t^{(r-q)} \int_{\Omega} b(x)|u|^{r}.
\]
Since \( H_u(t) - (q-\delta) \|u\|_{\wsq}^{q} = -G_u(t) \), it follows that \( H_u(t) < 0 \) for sufficiently small \( t \) and that \( H_u(t) \to \infty \) as \( t \to \infty \). Consequently, there exists a unique \( t_* > 0 \) satisfying \( H_u(t_*) = 0 \). In particular, for any fixed \( \lambda \in (0, \lambda_0) \), we obtain  
\[
t_* = \left(\frac{(p-\delta)\|u\|_{\wop}^{p}}{(r-\delta)\lambda \int_{\Omega} b(x)|u|^{r}}\right)^{\frac{1}{r-p}} > 0.
\]
Thus, there exists a unique \( t_{\max} > t_* > 0 \) such that \( H_u(t_{\max}) = (q-\delta) \|u\|_{\wsq}^{q} \). Furthermore, the function \( m_u(t) \) is increasing for \( t \in (0, t_{\max}) \) and decreasing for \( t \in (t_{\max}, \infty) \).
Consequently, 
\begin{equation*}
    \begin{aligned}
 (p-\delta) t_{max}^{p} \|u\|_{\wop}^{p} \leq    (p-\delta) t_{max}^{p} \|u\|_{\wop}^{p}&+(q-\delta) t_{max}^{q}\|u\|_{\wsq}^{q}
 \\
&=\lambda (r-\delta)t_{max}^{r} \int_{\Omega}  b(x)|u|^{r}
  \leq \lambda (r-\delta) t_{max}^{r} \bnorm S_{rp}^{-r/p}\|u\|_{\wop}^{r}. 
 \end{aligned}
\end{equation*}
Define $T_{0} = \frac{1}{ \|u\|_{\wop}} \left(\frac{p-\delta}{r-\delta} \frac{S_{rp}^{r/p}}{\lambda\bnorm}\right)^{\frac{1}{r-p}} \leq t_{\max},$ and
\begin{equation*}
    \begin{aligned}
        m_{u}(t_{\max}) &\geq  m_{u}(T_{0}) \geq T_{0}^{p-\delta}  \|u\|_{\wop}^{p} - T_{0}^{r-\delta} \lambda \bnorm  S_{p}^{\frac{-r}{p}}\|u\|_{\wop}^{r} \\
        &= \|u\|_{\wop}^{\delta} \left(\frac{r-p}{r-\delta}\right) \left(\frac{p-\delta}{r-\delta} \frac{S_{p}^{r/p}}{\bnorm}\right)^{\frac{p-\delta}{r-p}} \left(\frac{1}{\lambda}\right)^{\frac{p-\delta}{r-p}}.
    \end{aligned}
\end{equation*}
Therefore, if \( \lambda < \lambda_{0} \), then $\lambda\int_{\Omega} a(x)|u|^{\delta} \leq \lambda_{0} \anorm S_{p}^{\frac{-\delta}{p}} \|u\|_{\wop}^{\delta} \leq  m_{u}(t_{\max}).$ Thus, there exist \( t_{1} < t_{\max} \) and \( t_{2} > t_{\max} \) such that 
\[
m_{u}(t_{1}) = m_{u}(t_{2}) = \lambda\int_{\Omega} a(x)|u|^{\delta}.
\]
That is, \( t_{1}u, t_{2}u \in \nl \). Since \( m_{u} \) is increasing in \( (0, t_{\max}) \), it follows that \( m_{u}^{\prime}(t_{1}) > 0 \) and \( m_{u}^{\prime}(t_{2}) < 0 \), leading to \( t_{1}u \in \nl^{+} \) and \( t_{2}u \in \nl^{-} \). Now, 
\[
\gammau^{\prime}(t) = t^{\delta} \left(m_{u}(t) - \lambda\int_{\Omega}a(x)|u|^{\delta}\right),
\]
and utilizing the monotonicity properties of \( m_{u} \), we obtain 
\[
\gammau^{\prime }(t) < 0 \text{ } \text{for all } t \in [0, t_{1}),\quad
\gammau^{\prime }(t) > 0 \text{ } \text{for all } t \in (t_{1}, t_{2}].
\]
Thus, 
\[
J_{\lambda}(t_{1}u) = \min_{t \in [0, t_{2}]} J_{\lambda}(tu).
\]
Moreover,
\[
\gammau^{\prime}(t) > 0 \text{ } \text{for all } t \in [t_{1}, t_{2}),
\quad
\gammau^{\prime}(t) = 0 \text{ } \text{for } t = t_{2},
\quad
\gammau^{\prime}(t) < 0 \text{ } \text{for all } t \in [t_{2}, \infty).
\]
This implies
\[
J_{\lambda}(t_2 u) = \max_{t \geq t_{\max} } J_{\lambda} (t u).
\]
The case of $q>p $ can be handled by rewriting 
\begin{equation*}
    \begin{aligned}
          m_{u}^{\prime}(t)=t^{(p-\delta-1)}\left[(p-\delta)  \|u\|_{\wop}^{p}+(q-\delta) t^{(q-p)} \|u\|_{\wsq}^{q}-\lambda (r-\delta)t^{(r-p)} \int_{\Omega}  b(x)|u|^{r}\right]
    \end{aligned}
\end{equation*}
and following a similar arguments as above.
\\
\textbf{Case $(ii)$:} Suppose \( \int_{\Omega} a(x) |u|^{\delta} < 0 \) and \( \int_{\Omega} b(x) |u|^{r} > 0 \).  
We observe that \( m_{u}(t) \to -\infty \) as \( t \to \infty \), while \( m_{u}(t) > 0 \) for sufficiently small \( t \), and \( m_{u}^{\prime}(t) < 0 \) for large \( t \).  
By an argument similar to Case 1, there exists a unique \( t_0 > 0 \) such that \( m_{u} \) is increasing on \( (0,t_0) \) and decreasing on \( (t_0,\infty) \), with \( m_{u}^{\prime}(t_0) = 0 \).  
Since \( m_{u}(t_0) > 0 \) and \( \lambda \int_{\Omega} a(x) |u|^{\delta} < 0 \), there exists a unique \( t_1 > 0 \) such that  
\[
m_{u}(t_1) = \lambda \int_{\Omega} a(x) |u|^{\delta}, \quad m_{u}^{\prime}(t_1) < 0.
\]
This implies that \( t_1 u \in \nl^{-} \), meaning \( t_1 u \) is a local maximum.
\\
\textbf{Case $(iii)$:} Suppose \( \int_{\Omega} a(x) |u|^{\delta} > 0 \) and \( \int_{\Omega} b(x) |u|^{r} < 0 \).  
In this case, \( m_{u}^{\prime}(t) > 0 \) for all \( t > 0 \), meaning that \( m_{u} \) is an increasing function.  
Thus, there exists a unique \( t_1 > 0 \) such that  \[
m_{u}(t_1) = \lambda \int_{\Omega} a(x) |u|^{\delta}.
\]
Since \( \gamma_{t_1 u}^{\prime\prime}(1) > 0 \), we conclude that \( t_1 u \in \nl^{+} \), meaning \( t_1 u \) is a local minimum.
\\
\textbf{Case $(iv)$:} Suppose \( \int_{\Omega} a(x) |u|^{\delta} < 0 \) and \( \int_{\Omega} b(x) |u|^{r} < 0 \).  
In this case, we have \( \gammau(0) = 0 \) and \( \gammau^{\prime}(t) > 0 \) for all \( t > 0 \), implying that \( \gammau \) is strictly increasing and has no critical point. 
\end{proof}

\section{The Palais–Smale condition}\label{pscnd}
In this section, we derive several auxiliary results that will aid in establishing the Palais–Smale condition for the functional \( J_\lambda \) on the Nehari manifold. More generally, consider a Banach space \( X \) with a given norm and a functional \( I: X \to \mathbb{R} \) of class \( C^1 \). A sequence \( (u_n) \subset X \) is said to be a Palais–Smale sequence at level \( c \in \mathbb{R} \), abbreviated as \( (PS)_c \), if it satisfies the conditions \( I(u_n) \to c \) and \( I'(u_n) \to 0 \) as \( n \to \infty \). The Palais–Smale condition at level \( c \), or the \( (PS)_c \) condition, holds when every such sequence has a convergent subsequence. If this property is valid for all \( c \in \mathbb{R} \), we simply state that \( I \) satisfies the Palais–Smale condition.

\begin{lemma}\label{thetalambdaplus}
   There exist constants \( C_2, C_3 > 0 \) such that 
   \begin{equation*}
  \theta_{\lambda}^{+}\leq \left\{  \begin{aligned}
& -\frac{(p-\delta)(r-p)}{p\delta r} C_2 \hspace{3.9cm}\;\;\;\;\;\; \text{ if } p>q,
 \\
 & \max\left\{-\frac{(p-\delta)(r-p)}{p\delta r},-\frac{(q-\delta)(r-q)}{q\delta r}\right\} C_3  \;\;\;\;\;\; \text{ if }q \leq p.
 \end{aligned}\right.
\end{equation*}    Moreover, \( \theta_{\lambda}^{+} < 0 \).
\end{lemma}

\begin{proof}
Let \( u_{0} \in \wo(\Omega) \) be such that \( \int_{\Omega} a(x) |u_0|^{\delta} > 0 \).  
Then from previous lemma, there exists \( t_{0} > 0 \) such that \( t_{0} u_{0} \in \nl^{+} \), i.e., \( \gamma_{t_0 u_0}^{\prime\prime}(1) > 0 \). Since \( t_0 u_0 \in \nl \) and \( \gamma_{t_0 u_0}^{\prime\prime}(1) > 0 \), we obtain
\begin{equation}
  \begin{aligned}
J_{\lambda}(t_0 u_0) &= \left(\frac{1}{p} - \frac{1}{r} \right) \|t_0 u_0\|_{W_0^{1,p}}^{p}  
+ \left(\frac{1}{q} - \frac{1}{r} \right) \|t_0 u_0\|_{W_0^{s,q}}^{q} - \lambda \left(\frac{1}{\delta} - \frac{1}{r} \right) \int_{\Omega} a(x) |t_0 u_0|^{\delta} \\
&\leq \left(\frac{1}{p} - \frac{1}{r} \right) \|t_0 u_0\|_{W_0^{1,p}}^{p}  
+ \left(\frac{1}{q} - \frac{1}{r} \right) \|t_0 u_0\|_{W_0^{s,q}}^{q} \\
&\hspace{3.9cm} - \left(\frac{1}{\delta} - \frac{1}{r} \right) \left[ \frac{r-p}{r-\delta} \|t_0 u_0\|_{W_0^{1,p}}^{p} + \frac{r-q}{r-\delta} \|t_0 u_0\|_{W_0^{s,q}}^{q} \right].
  \end{aligned}
\end{equation}
Thus, if \( p > q \), then
\begin{equation}
  \begin{aligned}
J_{\lambda}(t_0 u_0) &\leq \left(\frac{1}{p} - \frac{1}{r} - \frac{r-p}{r\delta} \right) \|t_0 u_0\|_{W_0^{1,p}}^{p} = -\frac{(p-\delta)(r-p)}{p\delta r}  \|t_0 u_0\|_{W_0^{1,p}}^{p} < 0.
  \end{aligned}
\end{equation}
Similarly, if \( q > p \), we obtain
\begin{equation}
  \begin{aligned}
J_{\lambda}(t_0 u_0) &\leq \max\left\{-\frac{(p-\delta)(r-p)}{p\delta r}, -\frac{(q-\delta)(r-q)}{q\delta r} \right\} \left( \|t_0 u_0\|_{W_0^{1,p}}^{p} + \|t_0 u_0\|_{W_0^{s,q}}^{q} \right) < 0.
  \end{aligned}
\end{equation}
Therefore, the result follows.
\end{proof}
\begin{lemma}
    \label{thetalambdaminus}
    There exists $c_1>0$ such that $\jl(u)\geq c_1$ for any $u\in\nl^{-}$ for $\lambda\in(0,\widetilde\lambda)$ small enough. In particular $\theta_{\lambda}^{-}>0.$
\end{lemma}
\begin{proof}
    Let $u\in\nl^{-}.$ Then $\gammau^{\prime}(1)=0$ and $\gammau^{\prpr}(1)<0.$ Thus we have 
    \begin{equation*}
        \begin{aligned}
            (p-\delta)\|u\|_{\wop}^{p}+(q-\delta)\|u\|_{\wsq}^{q}<\lambda(r-\delta)\int_{\Omega}b|u|^{r}.
        \end{aligned}
    \end{equation*}
    Consequently we get there exists $A(\lambda_{0})>0$ such that $ \|u\|_{\wo}>A>0$ holds for all $u\in \nl^-, \; \lambda\in (0,\lambda_0).$  Recall that
    \begin{equation*}
        \begin{aligned}
            J_{\lambda}( u) &= \left(\frac{1}{p} - \frac{1}{r} \right) \|u\|_{W_0^{1,p}}^{p}  
+ \left(\frac{1}{q} - \frac{1}{r} \right) \| u\|_{W_0^{s,q}}^{q} - \lambda \left(\frac{1}{\delta} - \frac{1}{r} \right) \int_{\Omega} a(x) | u|^{\delta}
\\
&\geq \|u\|_{W_0^{1,p}}^{\delta}\left[\left(\frac{1}{p} - \frac{1}{r} \right) A^{p-\delta}  - \lambda \left(\frac{1}{\delta} - \frac{1}{r} \right) \anorm S_{rp}^{-\delta/p}\right] .
        \end{aligned}
    \end{equation*}
    Now we can choose $\widetilde\lambda$ small enough and obtain the result.
\end{proof}

\begin{lemma}\label{projlem}
   Let \( \lambda \in (0, \lambda_0) \) and \( z \in \nl \). Then there exist \( \varepsilon > 0 \) and a differentiable function \( \xi: B(0, \varepsilon) \subset \wo(\Omega) \to \mathbb{R} \) such that \( \xi(0) = 1 \), \( \xi(w)(z-w) \in \nl \), and 
   \begin{equation}\label{xiderivative}
       \begin{aligned}
           \langle \xi^{\prime}(0),w \rangle = 
           \frac{p \mathfrak{A}_{p}(z,w) + q A_{q}(z,w,\mathbb{R}^{2N}) - \lambda \int_{\Omega} \left( \delta a(x) [z]^{\delta-1}w + r b(x) [z]^{r-1}w \right)}
           {(p-\delta) \|z\|_{W_0^{1,p}}^{p} + (q-\delta) \|z\|_{W_0^{s,q}}^{q} - \lambda (r-\delta) \int_{\Omega} b(x) |z|^{r} }.
       \end{aligned}
   \end{equation}
   for all $w\in\wo(\Omega).$
\end{lemma}

\begin{proof}
For \( z \in \nl \), define the function \( H_z: \mathbb{R} \times \wo(\Omega) \to \mathbb{R} \) by 
\begin{equation*}
    \begin{aligned}
        H_{z}(t,w) &= \langle \jl^{\prime}(t(z-w)), t(z-w) \rangle \\
        &= t^{p} \|z-w\|_{W_0^{1,p}}^{p} + t^{q} \|z-w\|_{W_0^{s,q}}^{q} - \lambda \int_{\Omega} \left( a(x) t^{\delta} |z-w|^{\delta} + b(x) t^{r} |z-w|^{r} \right).
    \end{aligned}
\end{equation*}
Since \( H_{z}(1,0) = \langle \jl^{\prime}(z), z \rangle = 0 \) and \( \nl^{0} = \varnothing \), we obtain 
\begin{equation*}
    \frac{\partial}{\partial t} H_{z}(1,0) = (p-\delta) \|z\|_{W_0^{1,p}}^{p} + (q-\delta) \|z\|_{W_0^{s,q}}^{q} - \lambda (r-\delta) \int_{\Omega} b(x) |z|^{r}\not=0.
\end{equation*}
By the implicit function theorem, there exist \( \varepsilon > 0 \) and a differentiable function \( \xi: B(0, \varepsilon) \subset \wo(\Omega) \to \mathbb{R} \) such that \( \xi(0) = 1 \), satisfying \eqref{xiderivative}, and fulfilling \( H_{z}(\xi(w), w) = 0 \) for \( w \in B(0, \varepsilon) \), which is equivalent to  
\begin{equation*}
    \|\xi(w)(z-w)\|_{W_0^{1,p}}^{p} + \|\xi(w)(z-w)\|_{W_0^{s,q}}^{q} 
    - \lambda \int_{\Omega} \left( a(x) |\xi(w)(z-w)|^{\delta} + b(x) |\xi(w)(z-w)|^{r} \right) = 0.
\end{equation*}
Thus, for all \( w \in B(0, \varepsilon) \), we have \( \xi(w)(z-w) \in \nl \).
\end{proof}
 \begin{proposition}\label{criticallevel}
   Let  $\lambda \in (0,\lambda_0)$. Then, there exists a sequence $\{u_k\} \subset \nl$ such that
   \begin{equation*}
       \jl(u_k) = \theta_\lambda + o_k(1), \quad \jl^{\prime}(u_k) = o_k(1).
   \end{equation*}
\end{proposition}

\begin{proof}
    Since $\jl$ is coercive and bounded below in $\nl$, the Ekeland variational principle guarantees the existence of a minimizing sequence $\{u_k\} \subset \nl$ satisfying
   \begin{equation}\label{ekelandseq}
       \jl(u_k) < \theta_\lambda + \frac{1}{k}, \quad \jl(u_k) < \jl(v) + \frac{1}{k} \|v - u_k\|_{\wo}, \quad \forall v \in \nl.
   \end{equation}
   Since $u_k \in \nl$, we have
   \begin{equation*}
       \jl(u_k) = \left(\frac{1}{p} - \frac{1}{r}\right) \|u_k\|_{\wop}^{p} + \left(\frac{1}{q} - \frac{1}{r}\right) \|u_k\|_{\wsq}^{q} - \lambda \left(\frac{1}{\delta} - \frac{1}{r}\right) \int_\Omega a(x) |u_k|^{\delta}.
   \end{equation*}
   By \eqref{ekelandseq} and Lemma \ref{thetalambdaplus},
   \begin{equation*}
    \begin{aligned}
        \jl(u_k) <\theta_\lambda+ \frac{1}{k}\leq\theta_\lambda^{+}+ \frac{1}{k}<0 \text{ for sufficiently large } k.
    \end{aligned}
\end{equation*} 
It follows that $u_k \not\equiv 0$ for sufficiently large $k$. Moreover, applying H\"older's inequality, we obtain uniform bounds on $\|u_k\|_{\wop}$ and $\|u_k\|_{\wsq}$ i.e.
\begin{equation*}
    \begin{aligned}
        \|u_k\|_{\wop}\leq \left(\frac{\lambda(r-\delta)p\anorm}{(r-p)\delta S_{rp}^{\frac{\delta}{p}}}\right)^{\frac{1}{p-\delta}} \text{ and }   \|u_k\|_{\wsq}\leq \left(\frac{\lambda(r-\delta)q\anorm}{(r-q)\delta S_{rq}^{\frac{\delta}{q}}}\right)^{\frac{1}{q-\delta}}.
    \end{aligned}
\end{equation*}
Observe that $-\theta_\lambda^{+} \leq \sup_{\nl^{+}}\lambda\left(\frac{1}{\delta}-\frac{1}{r}\right)\int_{\Omega}a(x)|u|^{\delta}$. In particular, for the sequence $u_k$, we obtain the inequalities:
\begin{equation*}
    \begin{aligned}
        \left(\frac{(-\theta_\lambda^{+})\delta r S_{rp}^{\frac{\delta}{p}}}{(r-\delta)\anorm\lambda}\right)^{\frac{1}{\delta}} \leq  \|u_k\|_{\wop} \quad \text{and} \quad \left(\frac{(-\theta_\lambda^{+})\delta r S_{rq}^{\frac{\delta}{q}}}{(r-\delta)\anorm\lambda}\right)^{\frac{1}{\delta}} \leq  \|u_k\|_{\wsq}.
    \end{aligned}
\end{equation*}
Next, we aim to prove that $\|\jl^{\prime}(u_k)\|\rightarrow 0$ as $k\to\infty$. Using
Lemma \ref{projlem}, for each $u_k$, there exist $\varepsilon_k>0$ small and differentiable functions $\xi_k: B(0, \varepsilon_k) \subset \wo \to \mathbb{R}_{+}$ with $\xi_k(0) = 1$ and $\xi_k(v)(u_k - v) \in \nl$ for all $v \in B(0, \varepsilon_k)$. Fix $k\in\mathbb{N}$ such that $u_{k}\not\equiv 0$ and $0<\rho<\varepsilon_k.$ Setting $v_\rho = \frac{\rho u}{\|u\|_{\wo}}$ for an arbitrary $u\in\wo$ and $h_\rho = \xi_k(v_\rho)(u_k - v_\rho)$, we deduce from \eqref{ekelandseq} that
   \begin{equation*}
       \jl(h_\rho) - \jl(u_k) \geq -\frac{1}{k} \|h_\rho - u_k\|_{\wo}.
   \end{equation*}
Applying Taylor's theorem around $u_k$, we obtain
   \begin{equation*}
       \langle \jl^{\prime}(u_k), h_\rho - u_k \rangle + o_k(\|h_\rho - u_k\|_{\wo}) \geq -\frac{1}{k} \|h_\rho - u_k\|_{\wo^{\prime}}.
   \end{equation*}
Substituting $h_\rho - u_k = -v_\rho + (\xi_k(v_\rho) - 1)(u_k - v_\rho)$ and simplifying using   $\lim_{\rho\rightarrow 0 }\frac{|\xi_{k}(v_\rho)-1|}{\rho}\leq \|\xi_{k}^{\prime}(0)\|_{\wo^\prime},$ $h_\rho\in\nl,$ and $\|h_\rho-u_k\|_{\wo}\leq \rho|\xi_{k}(v_\rho)|+|\xi_{k}(v_\rho)-1|\|u_{k}\|_{\wo} $, we find that
\begin{equation*}
    \begin{aligned}
        \left\langle\jl^\prime(u_k),\frac{u}{\|u\|_{\wo}}\right\rangle \leq \frac{(\xi_k(v_\rho)-1)}{\rho}\langle\jl^{\prime}(u_k)-\jl^{\prime}(h_\rho),(u_k-v_{\rho})\rangle&+\frac{1}{k\rho} \|h_{\rho}-u_{k}\|_{\wo} \\&+\frac{1}{\rho} o_{k}(\|h_{\rho}-u_{k}\|_{\wo}).
    \end{aligned}
\end{equation*}
Thus for fixed $k$ and letting $\rho\rightarrow0$  in the above inequality, we observe that the first and last terms of the RHS converges to $0.$ The middle term can be estimated from above and we obtain 
\begin{equation}\label{jlambda_k}
       \left\langle \jl^{\prime}(u_k), \frac{u}{\|u\|_{\wo}} \right\rangle \leq \frac{C}{k} (1 + \|\xi_k^{\prime}(0)\|_{\wo^\prime})
   \end{equation}
   for some $C > 0$ independent of $k$. Next we claim that $\|\xi_k^{\prime}(0)\|_{\wo^\prime}$ is bounded. From Lemma \ref{projlem}, we estimate
 \begin{equation*}
       \begin{aligned}
           \langle \xi_k^{\prime}(0),v\rangle \leq \frac{K(\lambda_{0})\|v\|_{\wo}}{(p-\delta) \|u_k\|_{\wop}^{p}+(q-\delta) \|u_k\|_{\wsq}^{q}-\lambda (r-\delta) \int_{\Omega}  b(x)|u_k|^{r}}
       \end{aligned}
   \end{equation*}
   where $K( \lambda_0) > 0$. Suppose, for contradiction, that the denominator approaches zero along a subsequence. Then, since $u_k \in \nl$, we obtain $E_\lambda(u_k) = o_k(1).$ Furthermore, we have the following bounds for the norms of \( u_k \):
\begin{equation*}
    \begin{aligned}
        \|u_k\|_{\wop} &\geq \left( \frac{(p-\delta)S_{rp}^{\frac{r}{p}}}{\lambda_0(r-\delta)\bnorm} \right)^{\frac{1}{r-p}} + o_{k}(1), \quad
        \|u_k\|_{\wsq} \geq \left( \frac{(q-\delta)S_{rq}^{\frac{r}{q}}}{\lambda_0(r-\delta)\bnorm} \right)^{\frac{1}{r-q}} + o_{k}(1).
    \end{aligned}
\end{equation*}
Hence, there exists a constant \( d > 0 \) such that \( \|u_k\|_{\wo} > d > 0 \) for sufficiently large \( k \). Recall that 
\begin{equation*}
    E_{\lambda}(u_k) \geq \frac{r-p}{r-\delta} \|u_k\|_{\wop}^{p} - \lambda \anorm S_{rp}^{-\frac{\delta}{p}} \|u_k\|_{\wop}^{\delta}.
\end{equation*}
This leads to \( E_\lambda(u_k) > 0 \) for large \( k \), which contradicts the fact that \( E_\lambda(u_k) = o_k(1) \). Thus, the claim is established.
Since $\|\xi'_k(0)\|_{\mathcal{W}_0'}$ is bounded, from $\eqref{jlambda_k}$ we conclude that $\jl'(u_k)=o_k(1).$
\end{proof}
\section{Multiplicity Result When $p\leq q$ and $r<\max\{\pt,q_{s}^{*}\}$ } \label{subcritical}
In this section, we consider the case \( 1 < p \leq q < r < \max\{\pt,q_{s}^{*}\} \) and establish the existence and multiplicity results.

\begin{lemma}\label{psseq}
    If $\{u_k\}\subset\wo$ satisfies 
    \begin{equation*}
        \jl(u_k)=c+o_{k}(1), \quad \jl^{\prime}(u_k)=o_k(1) \text{ in } \w^{\prime},
    \end{equation*}
    then $\{u_k\}$ admits a convergent subsequence in $\wo.$
\end{lemma}
\begin{proof}
    Observe that $\{u_k\}$ is bounded in $\wo.$ Thus there exists $u_\lambda\in\wo$ such that, up to a subsequence, $u_k \rightharpoonup u_\lambda$ in $\wo$, $u_k \to u_\lambda$ strongly in $L^{\gamma}$ for $1\leq \gamma<\max\{\pt,q_{s}^{*}\}$, and $u_k(x) \to u_\lambda(x)$ a.e. in $\Omega$. The condition $\langle \jl^{\prime}(u_k)-\jl^{\prime}(u_\lambda),u_k-u_\lambda\rangle \to 0$ as $k\to\infty$ yields 
     \begin{equation*}
        \begin{aligned}
            o_k(1)
            =&\langle \jl^{\prime}(u_k)-\jl^{\prime}(u_\lambda),u_k-u_\lambda\rangle
            \\=&\mathfrak{A}_{p}(u_k,u_k-u_\lambda)-\mathfrak{A}_{p}(u_\lambda,u_k-u_\lambda)+A_{q}(u_k,u_k-u_\lambda,\RR^{2N})-A_{q}(u_\lambda,u_k-u_\lambda,\RR^{2N})
            \\
            &-\lambda\int_{\Omega}\left[a([u_k]^{\delta-1}-[u_\lambda]^{\delta-1})(u_k-u_\lambda)+b([u_k]^{r-1}-[u_\lambda]^{r-1})(u_k-u_\lambda)\right].
        \end{aligned}
    \end{equation*}
    Using H\"older's inequality, we obtain 
    \begin{equation*}
        \begin{aligned}
            \int_{\Omega} a(x)[u_k]^{\delta-1} (u_k-u_\lambda) &\leq \anorm \|u_k\|_{L^{r}}^{\delta-1}\|u_k-u_\lambda\|_{L^{r}} \to 0, \\
            \int_{\Omega} b(x)[u_k]^{r-1} (u_k-u_\lambda) &\leq \bnorm \|u_k\|_{L^{r}}^{r-1}\|u_k-u_\lambda\|_{L^{r}} \to 0.
        \end{aligned}
    \end{equation*}
Consider first the case $2\leq p\leq q$. Using the inequality 
    \begin{equation*}
        \begin{aligned}
            |a-b|^{l}\leq 2^{l-2}([a]^{l-1}-[b]^{l-1})(a-b) \text{ for }l\geq 2 \text{ and } a,b\in\RR^{N}.
        \end{aligned}
    \end{equation*}
    we obtain 
    \begin{equation*}
        \|u_k-u_\lambda\|_{\wo} \leq \langle \jl^{\prime}(u_k)-\jl^{\prime}(u_\lambda),u_k-u_\lambda\rangle =o_k(1),
    \end{equation*}
    implying strong convergence in $\wo$. Now, for $1<p\leq q<2$, using  the inequality,
     \begin{equation*}
       \begin{aligned}
           |a-b|^{l}\leq C_{l}(([a]^{l-1}-[b]^{l-1})(a-b))^{l/2}(|a|^{l}+|b|^{l})^{\frac{2-l}{2}}\text{ for }1<l< 2 \text{ and } a,b\in\RR^{N},
       \end{aligned}
   \end{equation*}
    we deduce 
    \begin{equation*}
        \|u_k-u_\lambda\|_{\wop}^{p} \leq C(\mathfrak{A}_{p}(u_k,u_k-u_\lambda)-\mathfrak{A}_{p}(u_\lambda,u_k-u_\lambda))^{p/2} \big(\|u_k\|_{\wop}^{p} + \|u_\lambda\|_{\wop}^{p}\big)^{(2-p)/2}.
    \end{equation*}
    Since $\{u_k\}$ is bounded in $\wo$, it follows that 
    \begin{equation*}
        \|u_k-u_\lambda\|_{\wop}^{2} \leq C(\mathfrak{A}_{p}(u_k,u_k-u_\lambda)-\mathfrak{A}_{p}(u_\lambda,u_k-u_\lambda)).
    \end{equation*}
    Similarly, 
    \begin{equation*}
        \|u_k-u_\lambda\|_{\wsq}^{2} \leq C (A_{q}(u_k,u_k-u_\lambda,\RR^{2N})-A_{q}(u_\lambda,u_k-u_\lambda,\RR^{2N})).
    \end{equation*}
    Combining these, we conclude $\|u_k-u_\lambda\|_{\wo} \to 0$ as $k\to\infty$. The argument extends similarly to the case $1<p<2<q$, completing the proof.
\end{proof}
Now we prove the existence of two solution for the subcritical case.
\\[1mm]
\textit{\underline{\textbf{Proof of Theorem \ref{subcrmul}} }}
\\[3mm]
    By Proposition \ref{criticallevel}, there exist minimizing sequences \( \{u_k\} \subset \nl^{+} \) and \( \{v_k\} \subset \nl^{-} \) with respect to $\jl.$ In the view of Lemma \ref{thetalambdaplus} and Lemma \ref{thetalambdaminus} we can apply Lemma \ref{psseq} in $\nl^\pm$ and deduce the existence of \( u_\lambda, v_\lambda \in \wo \) such that \( u_k \to u_\lambda \) and \( v_k \to v_\lambda \) strongly in \( \wo \) for all \( \lambda \in (0, \lambda_0) \). Consequently, \( u_\lambda \) and \( v_\lambda \) are weak solutions of \eqref{generaleqn}. 
    \\
From Lemma \ref{thetalambdaplus}, it follows that \( u_\lambda \not\equiv 0 \), ensuring \( u_\lambda \in \nl \). Furthermore, Lemma \ref{nlambdazero} implies \( u_\lambda \in \nl^{+} \) with \( \jl(u_\lambda) = \theta_{\lambda}^{+}<0 \).
Since \( \jl(v_\lambda) = \theta_{\lambda}^{-}>0 \) and \( \nl^{+} \cap \nl^{-} = \varnothing \), we have \( v_\lambda \in \nl^{-} \). And consequently, the solutions \( u_\lambda \) and \( v_\lambda \) are distinct.
\\
Next, we establish the non-negativity of \( u_\lambda \). If \( u_\lambda \geq 0 \), it is already a nonnegative solution of \eqref{generaleqn} and a minimizer for \( \jl \) in \( \nl^{+} \). Otherwise, by Lemma \ref{fiberingmapanalysis}, there exists a unique \( t_1 > 0 \) such that \( t_1 u_\lambda \in \nl^{+} \). Observing that 
\[
M_{|u_{\lambda}|}(1) \leq M_{u_{\lambda}}(1) = \lambda \int_\Omega a|u_\lambda|^{\delta} \leq M_{|u_{\lambda}|}(t_1) \leq M_{u_{\lambda}}(t_1),
\]
and noting that \( M_{u_{\lambda}}^{\prime}(1) > 0 \) (since \( u_\lambda \in \nl^{+}\)), we infer \( t_1 \geq 1 \). Thus,
\[
\theta_{\lambda}^{+} \leq \gamma_{|u_{\lambda}|}(t_1) \leq \gamma_{u_{\lambda}}(1) = \theta_{\lambda}^{+}.
\]
This implies \( \jl(t_1 |u_{\lambda}|) = \gamma_{|u_{\lambda}|}(t_1) = \theta_{\lambda}^{+} \), and \( t_1|u_\lambda| \in \nl^{+} \). Therefore, \( t_1|u_\lambda| \) is a nonnegative solution of \eqref{generaleqn} in \( \nl^{+} \). A similar argument ensures that \( v_\lambda \) is also a nonnegative solution. \hfill\qed

\section{Nonexistence Result }\label{nonexistence}
Let \( sq < p  < \infty \) be fixed throughout this section. We begin by establishing a regularity result.
\begin{theorem}\label{regularity}
   Let  \( u \in \wo(\Omega) \) be a non-negative subsolution to the problem 
   \begin{equation}
    \left\{ \begin{aligned}
       & \fp u + \fqs u = |u|^{p-2}u + |u|^{q-2}u \quad \text{in } \Omega, \\
       & u = 0 \quad \text{in } \Omega^{c}.
    \end{aligned} \right.
   \end{equation}
   Then, \( \|u\|_{L^{\infty}(\Omega)} \) is bounded, depending only on \( N, p, q, s \), \( \|u\|_{L^{p}(\Omega)} \), and \( \|u\|_{L^{q}(\Omega)} \). Furthermore, \( u \in C^{1,\gamma}_{0}(\overline{\Omega}) \) for some \( \gamma \in (0,1) \).
\end{theorem}
\begin{sketch}
   Define \( f(s) = |s|^{p-1}s + |s|^{q-1}s \). It follows that \( f(s) \leq C(s^{-\delta} + s^{l}) \) for some \( 0 < \delta < 1 \) and \( \max\{p,q\} \leq l \leq \min\{\pt, \qst\} \). By adapting the proof of \cite[Theorem 7.1]{DGJ24}, we deduce that \( u \in L^{\infty}(\Omega) \). Subsequently, applying \cite[Theorem 4 and 5]{FM22} and \cite[Theorem 1.1]{AC23}, we conclude that \( u \in C^{1,\gamma}_{0}(\overline{\Omega}) \) for some \( \gamma \in (0,1) \).
\end{sketch}
To establish a nonexistence result for \eqref{generaleqn}, we first derive a nonexistence result for a generalized eigenvalue problem.
 \subsection{Nonexistence Result For Generalized Eigenvalue Problem}
We begin by considering the equation
\begin{equation*}
   (GEV,\alpha,\beta) \left\{ \begin{aligned}
       & \fp u + \fqs u = \alpha|u|^{p-2}u + \beta|u|^{q-2}u \quad \text{in } \Omega, \\
       & u = 0 \quad \text{in } \Omega^{c}.
    \end{aligned} \right.
\end{equation*}
First we consider the case $q<p.$ Define
\begin{equation*}
    \lambda^{*}(s) = \sup\left\{\lambda \in \RR : (GEV,\lambda,\lambda+s) \text{ has a positive solution}\right\}.
\end{equation*}
If no such \(\lambda\) exists for a fixed \(s \in \RR\), we set \(\lambda^{*}(s) = -\infty\). Our goal is to show that \(\lambda^{*}(s)\) is bounded independently of \(s\). To this end, we recall two key inequalities. From \cite[Proposition 8 and (9.5)]{BT15}, there exists \(\rho > 0\) such that
\begin{equation}\label{btapp}
    [\nabla u]^{t-1} \nabla \left(\frac{\varphi^{m}}{u^{t-1} + u^{m-1}}\right) \leq \frac{|\nabla (\varphi^{\frac{m}{t}})|^{t}}{\rho} \text{ if } 1<t<m<\infty,
\end{equation}
for any differentiable functions \(u > 0\) and \(\varphi \geq 0\). Additionally, by \cite[Remark 2.6]{GGM22}, we have
\begin{equation}\label{ggmeqn}
    [u(x) - u(y)]^{m-1} \left(\frac{\varphi^{m}}{u^{t-1} + u^{m-1}}(x) - \frac{\varphi^{m}}{u^{t-1} + u^{m-1}}(y)\right) \leq |\varphi(x) - \varphi(y)|^{m}.
\end{equation}
Using these inequalities, we now proceed to establish the boundedness of \(\lambda^{*}\).
\begin{lemma}\label{lastarbdd}
    For \( p < q \) and any fixed \( s \in \RR \), the value \( \lambda^{*}(s) \) is bounded.
\end{lemma}

\begin{proof}
    Fix \( s \in \RR \), and let \( u \in \wo(\Omega) \) be a positive solution of \( (GEV, \lambda, \lambda + s) \). By Theorem \ref{regularity} and \cite[Theorem 1.2]{AC23}, we have \( u \in \text{int}\, C_{0}^{1}(\overline{\Omega})_{+} \). Choose a test function \( \varphi \in \text{int}\, C_{0}^{1}(\overline{\Omega})_{+} \), and define \( \chi = \frac{\varphi^{q}}{u^{p-1} + u^{q-1}} \in \wo \). Substituting \( \chi \) into \( (GEV, \lambda, \lambda + s) \) and applying inequalities \eqref{btapp} and \eqref{ggmeqn}, we obtain
    \begin{equation*}
        \begin{aligned}
            \lambda \int_{\Omega} \varphi^{q} + s \int_{\Omega} \frac{u^{q-1} \varphi^q}{u^{q-1} + u^{p-1}} \leq \frac{1}{\rho} \int_{\Omega} |\nabla (\varphi^{\frac{q}{p}})|^{p} + \underset{\RR^{2N}}{\int\int} |\varphi(x) - \varphi(y)|^{q} \frac{dx\, dy}{|x-y|^{N+sq}}.
        \end{aligned}
    \end{equation*}
    Observe that \( \int_{\Omega} \frac{u^{q-1} \varphi^{q}}{u^{q-1} + u^{p-1}} \leq  \int_{\Omega} \varphi^{q} \). Thus, we derive
    \begin{equation*}
        \begin{aligned}
            \lambda \int_{\Omega} \varphi^{q} +  \min\left\{0, s\int_{\Omega} \varphi^{q}\right\} \leq \frac{1}{\rho} \int_{\Omega} |\nabla (\varphi^{\frac{q}{p}})|^{p} + \underset{\RR^{2N}}{\int\int} |\varphi(x) - \varphi(y)|^{q} \frac{dx\, dy}{|x-y|^{N+sq}}.
        \end{aligned}
    \end{equation*}
    Since \( \int_{\Omega} \varphi^{q} \), \( \int_{\Omega} |\nabla (\varphi^{\frac{q}{p}})|^{p} \), \( \|\varphi\|_{\wsq} \), and \( s \) are independent of \( \lambda \) and \( u \), it follows that if \( (GEV, \lambda, \lambda + s) \) admits a positive solution, \( \lambda \) must be bounded from above.
\end{proof}
Define \( f_{\alpha,\beta}(t) = \alpha [t]^{p-1} + \beta [t]^{q-1} \). For any two functions \( v, w \in L^{\infty}(\Omega) \) satisfying \( v \leq w \), we introduce the truncation
\begin{equation*}
    \widetilde{f}_{\alpha,\beta}^{[v,w]}(x,t) = \left\{
    \begin{aligned}
        & f_{\alpha,\beta}(w(x)) \quad \text{if } t \geq w(x), \\
        & f_{\alpha,\beta}(t) \quad \text{if } v(x) \leq t \leq w(x), \\
        & f_{\alpha,\beta}(v(x)) \quad \text{if } t \leq v(x).
    \end{aligned}\right.
\end{equation*}
The corresponding energy functional is defined as
\begin{equation*}
    E_{\alpha,\beta}^{[v,w]}(u) = \frac{1}{p}\|u\|_{\wop}^{p} + \frac{1}{q}\|u\|_{\wsq}^{q} - \int_{\Omega} \int_{0}^{u(x)} \widetilde{f}_{\alpha,\beta}^{[v,w]}(x,t) \, dt \, dx.
\end{equation*}

\begin{lemma}
    Let \( \alpha > \lambda_{1p} \) and $p<q$, where \( \lambda_{1p} \) is the first eigenvalue of \( \fp \). Moreover, let \( w \in \operatorname{int}\, C_{0}^{1}(\overline{\Omega})_{+} \) be a positive supersolution of \( (GEV, \alpha, \beta) \). Then,  \( \min_{\wo} E_{\alpha,\beta}^{[0,w]} < 0 \) and hence $(GEV, \alpha, \beta)$ has a positive solution in $\operatorname{int}(C_{0}^{1}(\overline{\Omega})_{+}).$
\end{lemma}

\begin{proof}
    The functional \( E_{\alpha,\beta}^{[v,w]} \) is coercive, weakly lower semicontinuous, and bounded below, ensuring the existence of a global minimum. Since \( w \) and \( \phi_{p} \) (the first eigenfunction of \( \fp \)) belong to \( \text{int}\, C_{0}^{1}(\overline{\Omega})_{+} \), we can choose \( t > 0 \) sufficiently small such that \( t\phi_{p} \leq w \). Consequently, \( \widetilde{f}_{\alpha,\beta}^{[0,w]}(x, t\phi_{p}) = f_{\alpha,\beta}(t\phi_{p}) \). Evaluating the energy functional at \( t\phi_{p} \), we obtain
    \begin{equation}\label{ealphabeta}
        \begin{aligned}
            E_{\alpha,\beta}^{[0,w]}(t\phi_{p}) &= \frac{t^{p}}{p}\|\phi_{p}\|_{\wop}^{p} + \frac{t^q}{q}\|\phi_{p}\|_{\wsq}^{q} - \int_{\Omega} \left( \frac{\alpha}{p} |t|^{p}|\phi_{p}|^{p} + \frac{\beta}{q} |t|^{q}|\phi_{p}|^{q} \right) \\
            &= \frac{t^p}{p}(\lambda_{1,p} - \alpha)\|\phi_{p}\|_{L^{p}} + \frac{t^{q}}{q}\left(\|\phi_{p}\|_{\wsq} - \beta \|\phi_{p}\|_{L^{q}}\right) \\
            &\leq \frac{t^p}{p}(\lambda_{1,p} - \alpha)\|\phi_{p}\|_{L^{p}} + \frac{t^{q}}{q}\|\phi_{p}\|_{\wsq}.
        \end{aligned}
    \end{equation}
    Since \( \alpha > \lambda_{1,p} \) and \( p < q \), for sufficiently small \( t \), the term \( \frac{t^p}{p}(\lambda_{1,p} - \alpha)\|\phi_{p}\|_{L^{p}} \) dominates, making \( E_{\alpha,\beta}^{[0,w]}(t\phi_{p}) < 0 \). This completes the proof.
\end{proof}

\begin{lemma}\label{unibdoflmbdstr}
    Let us consider the case $q<p$ and $\lambda^{*}(s)>\lambda_{1p}.$ Then we have the following properties in $\RR$:
    \begin{itemize}
        \item[(a)]  $\lambda^{*}(s)$ is nonincreasing; 
         \item[(b)]  $\lambda^{*}(s)+s$ is nondecreasing;
    \end{itemize}
\end{lemma}
\begin{proof}
\textit{Part (a)} We shall show that if both $\lambda^{*}(s)$ and $\lambda^{*}(s^{\prime})$ are larger than $\lambda_{1p},$ then $\lambda^{*}(s^{\prime})\leq \lambda^{*}(s)$ when $s<s^{\prime}.$  Fix any $\varepsilon>0$ such that $\lambda^{*}(s^{\prime})-\varepsilon>\lambda_{1p}.$  By the definition of $\lambda^{*}$ there exists $\mu$ such that $\lambda^{*}(s^{\prime})-\varepsilon<\mu<\lambda^{*}(s^{\prime})$ such that $(GEV,\mu,\mu+s^{\prime})$ has a positive solution $w_\mu \in int\;C_{0}^{1}(\overline{\Omega})_{+}.$ Since $s<s^{\prime}$ we get that $w_{\mu}$ is positive supersolution of $(GEV,\mu,\mu+s).$ By the previous lemma, $(GEV,\mu,\mu+s)$ has a positive solution. Then by definition $\lambda^{*}(s^{\prime})\leq \lambda^{*}(s).$
\\
\textit{Part (b)} follows from a similar argument as \textit{Part (a)}.(See \cite[Proposition 3]{BT15} for details.) 
\end{proof}
Next we consider the case $p<q$ and thus define 
\begin{equation*}
    \lambda_{1}^{*}(s) = \sup\left\{\lambda \in \RR : (GEV,\lambda+s,\lambda) \text{ has a positive solution}\right\}.
\end{equation*}
If no such \(\lambda\) exists for a fixed \(s \in \RR\), we set \(\lambda_{1}^{*}(s) = -\infty\).
Similar to \eqref{lastarbdd} we can prove the following lemma
\begin{lemma}
    For \( q < p \) and any fixed \( s \in \RR \), the value \( \lambda_{1}^{*}(s) \) is bounded.
\end{lemma}
In the case \( q < p \), one might naturally consider using the first eigenfunction \( \phi_q \) of \( \fqs \) as a test function for the energy functional \( E_{\alpha,\beta}^{[v,w]} \). However, due to the regularity properties of the fractional \( q \)-Laplacian, we have \( \phi_q \in C^s(\overline{\Omega}) \), which implies that the norm \( \|\phi_q\|_{\wop} \) may not be well-defined. Consequently, a more detailed analysis is required to address this case.
\begin{lemma}
    Let \( \beta > \lambda_{1q} \) and $q<p$, where \( \lambda_{1,q} \) is the first eigenvalue of \( \fqs \). Moreover, let \( w \in \operatorname{int}\, C_{0}^{1}(\overline{\Omega})_{+} \) be a positive supersolution of \( (GEV, \alpha, \beta) \). Then,  \( \min_{\wo} E_{\alpha,\beta}^{[0,w]} < 0 \) and hence $(GEV, \alpha, \beta)$ has a positive solution in $\operatorname{int}(C_{0}^{1}(\overline{\Omega})_{+}).$
\end{lemma}
\begin{proof}
    The functional \( E_{\alpha,\beta}^{[v,w]} \) is coercive, weakly lower semicontinuous, and bounded below, ensuring the existence of a global minimum. By \cite{FSV15} there exists a sequence $\{\phi_n\}\subset C_c^{\infty}(\Omega)$ such that $\phi_n\rightarrow \phi_q$ in $\wsq$ as $n\rightarrow\infty.$ Thus $|\phi_n|\rightarrow \phi_q$ in $\wsq$ as $n\rightarrow\infty.$ Since \( \beta > \lambda_{1,q} \), we have $\|\phi_{q}\|_{\wsq} - \beta \|\phi_{q}\|_{L^{q}}<0.$ Therefor there exists $n_0>0$ such that  $\|\phi_{n}\|_{\wsq} - \beta \|\phi_{n}\|_{L^{q}}<0$ for all $n\geq n_0.$ Since \( w \in \operatorname{int}\, C_{0}^{1}(\overline{\Omega})_{+} \) we can choose $t>0$ small enough such that  $0<t |\phi_{n_{0}}|\leq w$ in $\Omega.$ Now similar to \eqref{ealphabeta} we can conclude that $E_{\alpha,\beta}^{[0,w]}(t|\phi_{n_{0}}|)<0$ and complete the proof. 
    \end{proof}
    Now analogous to Lemma \ref{unibdoflmbdstr}, we can prove the following lemma.
    \begin{lemma}\label{unibdoflmbdstr1}
    Let us consider the case $q<p$ and $\lambda_{1}^{*}(s)>\lambda_{1q}.$ Then we have the following properties in $\RR$:
    \begin{itemize}
        \item[(a)]  $\lambda_{1}^{*}(s)$ is nonincreasing; 
         \item[(b)]  $\lambda_{1}^{*}(s)+s$ is nondecreasing;
    \end{itemize}
\end{lemma}

Using Part(a) and (b) of Lemma \ref{unibdoflmbdstr} and \ref{unibdoflmbdstr1}  along with the definition of $\lambda^*(s)$ we have the following result.
\begin{corollary}\label{nonex}
    For  a given $c_1,c_2>0$ the problem $(GEV,\mu c_1,\mu c_2)$ does not admit a positive solution if $\mu>>1.$
\end{corollary}


\subsection{Nonexistence result}
We now proceed to establish the nonexistence result for the subcritical problem with nonnegative weight function. We adopt the arguments of \cite{MMP19} in the mixed local nonlocal case. 
\\[1mm]
\textit{\underline{\textbf{Proof of Theorem \ref{nonexsubcr}} }}
\\[3mm]
    Assume, to the contrary, that \eqref{generaleqn} has a positive solution for any $\lambda > 0$. Fix such a $\lambda$ and let $u_\lambda \in \operatorname{int} C_0^1(\overline{\Omega})_+$ be a corresponding positive solution. Now there exists $c_1>0$ and $c_2>0$ such that 
    \begin{equation*}
        c_1 t^{p-1} + c_2 t^{q-1} \leq \alpha t^{\delta-1}+ \beta t^{r-1} \text{ for all } t\geq 0.
    \end{equation*}
    For that $c_1$ and $c_2$ we define the function
    \begin{equation}
        g_{\lambda}(x,t) = \begin{cases}
            \lambda c_1 t_+^{p-1} + \lambda c_2 t_+^{q-1}, & \text{if } t \leq u_{\lambda}, \\
            \lambda c_1 u_{\lambda}^{p-1} + \lambda c_2 u_{\lambda}^{q-1}, & \text{if } t > u_{\lambda}.
        \end{cases}
    \end{equation}
    Consider the associated energy functional
    \begin{equation*}
        \widetilde{E}_{\lambda}(u) = \frac{1}{p}\|u\|_{W_0^{s,p}}^p + \frac{1}{q}\|u\|_{W_0^{s,q}}^q - \int_{\Omega} \int_{0}^{u(x)} g_{\lambda}(x,t) \, dt \, dx.
    \end{equation*}
    Since $\widetilde{E}_{\lambda}$ attains a global minimum, let $\bar{u}_{\lambda}$ be a minimizer. Then, it satisfies the weak formulation
    \begin{equation}
        \mathfrak{A}_{p}(\bar{u}_{\lambda}, \varphi) + A_{q}(\bar{u}_{\lambda}, \varphi, \mathbb{R}^{2N}) = \int_{\Omega} g_{\lambda}(x, \bar{u}_{\lambda}(x)) \varphi(x) \, dx, \quad \forall \varphi \in \wo(\Omega).
    \end{equation}
    Choosing test functions $\varphi = (\bar{u}_{\lambda})_-$ and $\varphi = (\bar{u}_{\lambda} - u_{\lambda})_+$ yields $0 \leq \bar{u}_{\lambda} \leq u_{\lambda}$. Since $u_\lambda, \phi_p \in \operatorname{int} C_0^1(\overline{\Omega})_+$, there exists $t > 0$ small enough such that $t \phi_p \leq u_\lambda$. A computation analogous to \eqref{ealphabeta} shows
    \begin{equation*}
        \widetilde{E}_{\lambda}(t \phi_p) \leq \frac{t^p}{p}(\lambda_{1,p} - \lambda c_1) \|\phi_p\|_{L^p} + \frac{t^q}{q} \|\phi_p\|_{W_0^{s,q}} < 0 \quad \text{if } \lambda > \frac{\lambda_{1,p}}{c_1}.
    \end{equation*}
    Since $\bar{u}_\lambda$ is the global minimizer, it follows that $\bar{u}_\lambda$ satisfies
    \begin{equation*}
        \fp \bar{u}_\lambda + \fqs \bar{u}_\lambda = \lambda c_1 \bar{u}_\lambda^{p-1} + \lambda c_2 \bar{u}_\lambda^{q-1},
    \end{equation*}
    for any $\lambda$ large enough, contradicting Corollary \ref{nonex}. Thus, the claim holds. \hfill\qed
\begin{remark}  
A similar nonexistence result can be established for the case \( p = q \) by utilizing the isolation of the first eigenvalue of \( \fp + \fps \). In particular, under the same conditions as in Theorem \ref{nonexistence}, one can show the existence of a threshold \( \Lambda_* > 0 \) such that for all \( \lambda > \Lambda_* \), the problem \eqref{generaleqn} admits only the trivial solution when \( r < \max\{\pt, p_{s}^{*}\} = \pt \).  
\end{remark}  
This remark follows by contradiction using the isolation of the first eigenvalue \( \lambda_1 \) (see \cite[Theorem 1.1]{DFR19}). Suppose there exists a sequence \( \lambda_n \to \infty \) such that for each \( n \), problem \eqref{generaleqn} has a solution \( u_n \). Then, there exists \( \tilde{\lambda}_0 > 0 \) such that the following inequality holds:  
\begin{equation*}  
\lambda (\alpha t^{\delta-1} + \beta t^{r-1}) \geq (\lambda_1 + \varepsilon) t^{p-1}, \quad \text{for all } t > 0, \text{ and for any } 0 < \varepsilon < 1, \text{ with } \lambda > \tilde{\lambda}_0.  
\end{equation*}  
For \( \lambda_n > \tilde{\lambda}_0 \), the function \( u_n \) serves as a supersolution to the problem:  
\begin{equation}  
\label{subsupeqn}  
\begin{aligned}  
\fp u + \fps u &= (\lambda_1 + \varepsilon) u^{p-1}, \quad u > 0 \text{ in } \Omega, \quad u|_{\Omega^c} = 0,  
\end{aligned}  
\end{equation}  
for any \( 0 < \varepsilon < 1 \). By Hopf’s lemma, it follows that \( u_n \geq C(\tilde{\lambda}_0) \dis(x) \), where the constant \( C(\tilde{\lambda}_0) \) is independent of \( n \).  
Choosing \( \mu \) small enough so that \( \mu < \lambda_1 + \varepsilon \), we obtain \( \mu \phi_1 < C(\tilde{\lambda}_0) \dis(x) \leq u_n \), where \( \phi_1 \) is the normalized positive eigenfunction corresponding to \( \lambda_1 \). Moreover, \( \mu \phi_1 \) acts as a subsolution to \eqref{subsupeqn}. Applying the monotone iteration method then yields a solution to \eqref{subsupeqn} for any \( \varepsilon \in (0,1) \), which contradicts the fact that \( \lambda_1 \) is an isolated eigenvalue.

\section{Multiplicity Result When $q\leq p$ and $r=p_{*}$}\label{critical}
In this section, we assume $r = p_{*}$ and $q \leq  p.$ This implies that $\wo \equiv \wop(\Omega)$ and $\qst<\pt.$ Furthermore, we consider the case where $b(x) \equiv 1$ and $a(x)$ is a continuous function satisfying $\inf_{B_{r_0}(x_0)} a(x) = m_{a} > 0$ for some $r_0 > 0$. By a argument similar to \cite[Lemma 2.2]{DFV24} we get the following result.
\begin{lemma}
    Let \(\lambda > 0\) and let \(\{u_n\} \subset \wopo\) be a bounded $(PS)_c$ sequence  with \(c \in \mathbb{R}\).  Then, up to a subsequence, \(\nabla u_n (x) \to \nabla u(x)\) a.e. in \(\Omega\) as \(n \to \infty\).
\end{lemma}
\begin{theorem}\label{cdeltathm}
    Let $\lambda \in (0, \lambda_0)$, and suppose that $\{u_k\} \subset \nl$ is a $(PS)_c$ sequence for $\jl$, with $u_k$ converging weakly to $u$ in $\wop$. Then, $\jl^\prime(u)=0$ and  there exists a positive constant $C_{\delta}$, depending on $p, N, S_p, |\Omega|$, and $\delta$, such that
    \begin{equation}
        \jl(u) \geq -C_{\delta} \lambda^{\frac{p}{p - \delta}},
    \end{equation}
    where
    \begin{equation}\label{cdelta}
        C_{\delta} = \left( \frac{1}{\delta} - \frac{1}{p_{*}} \right) \left[ \left( \frac{p}{\delta} \left( \frac{1}{p} - \frac{1}{p_{*}} \right) \left( \frac{1}{\delta} - \frac{1}{p_{*}} \right)^{-1} \right)^{-\frac{\delta}{p}} \|a\|_{L^{\infty}} S_p^{-\frac{\delta}{p}} |\Omega|^{\frac{p_{*} - \delta}{p_{*}}} \right]^{\frac{p}{p - \delta}}.
    \end{equation}
\end{theorem}
\begin{proof}
    Since $u_k$ converges weakly to $u$ in $\wopo$, we conclude that $u_k \to u$ strongly in $L^{\gamma}$ for $1 \leq \gamma < p_{*}$, and pointwise almost everywhere in $\Omega$. By \cite[Lemma 2.2]{CMM18}, $\fqs$ is weak-weak continuous, implying that
    \begin{equation*}
        \lim_{k\to\infty} A_{q}(u_k, \varphi, \mathbb{R}^{2N}) = A_{q}(u, \varphi, \mathbb{R}^{2N}) \quad \text{for any } \varphi \in \wopo.
    \end{equation*}
    Since $u_k \rightharpoonup u$ weakly in $\wopo$, it follows that $[\nabla u_k]^{p-1}$ is bounded, and consequently,  $[\nabla u_k]^{p-1} \rightharpoonup (f_1,f_2\cdots f_n)$ in $(L^{p'})^n$. Since $\nabla u_k \rightarrow \nabla u$ a.e we conclude that $[\nabla u_k]^{p-1} \rightharpoonup [\nabla u]^{p-1} $ in $(L^{p'})^n.$ This implies, given a $\varphi \in \wopo$, we obtain
    \begin{equation*}
        \lim_{k \to \infty} \mathfrak{A}_{p}(u_k, \varphi) = \mathfrak{A}_{p}(u, \varphi).
    \end{equation*}
Moreover, using the weak convergence $u_k \rightharpoonup u$ in $\wopo$, we deduce that
    \begin{equation*}
         [u_k]^{\delta-1} \rightharpoonup [u]^{\delta-1} \text{ weakly in } L^{\delta'} \quad \text{and}
         \quad [u_k]^{\pt-1} \rightharpoonup [u]^{\pt-1} \text{ weakly in } L^{\pt^\prime}.
    \end{equation*}
    This yields the convergence
    \begin{equation*}
        \int_{\Omega} a(x)([u_k]^{\delta-1} - [u]^{\delta-1}) \varphi(x) \to 0, \quad \int_{\Omega} ([u_k]^{p_{*}-1} - [u]^{p_{*}-1}) \varphi(x) \to 0,
    \end{equation*}
    for any $\varphi \in \wopo \subset L^{\delta}(\Omega) \cap L^{p_{*}}(\Omega)$. Combining these results, we establish that
    \begin{equation*}
     \begin{aligned}
            \langle \jl'(u_k) - \jl'(u), \varphi \rangle &= \mathfrak{A}_{p}(u_k, \varphi) - \mathfrak{A}_{p}(u, \varphi) + A_{q}(u_k, \varphi, \mathbb{R}^{2N}) - A_{q}(u, \varphi, \mathbb{R}^{2N}) 
        \\&- \lambda \left( \int_{\Omega} a(x)([u_k]^{\delta-1} - [u]^{\delta-1}) \varphi(x) + ([u_k]^{p_{*}-1} - [u]^{p_{*}-1}) \varphi(x) \right) = o_k(1).
     \end{aligned}
    \end{equation*}
    Since $u_k$ is a $(PS)_c$ sequence, it follows that $\langle \jl'(u), \varphi \rangle = 0$, leading to
    \begin{equation}\label{jleqninnl}
        \jl(u) \geq \left( \frac{1}{p} - \frac{1}{p_{*}} \right) \|u\|_{\wop}^{p} - \lambda \left( \frac{1}{\delta} - \frac{1}{p_{*}} \right) \int_{\Omega} a(x) |u(x)|^{\delta} \,dx.
    \end{equation}
    Applying H\"older's inequality, Sobolev embeddings, and Young’s inequality, We derive
\begin{equation*}
    \begin{aligned}
        \lambda \int_{\Omega} a |u|^\delta &= \left( \frac{p}{\delta} \left( \frac{1}{p} - \frac{1}{\pt} \right) \left( \frac{1}{\delta} - \frac{1}{\pt} \right)^{-1} \right)^{\frac{\delta}{p}} \|u\|_{\wop}^{\delta}
        \\&
        \times \lambda \left( \frac{p}{\delta} \left( \frac{1}{p} - \frac{1}{\pt} \right) \left( \frac{1}{\delta} - \frac{1}{\pt} \right)^{-1} \right)^{\frac{-\delta}{p}} \|a\|_{L^{\infty}} S_{p}^{\frac{-\delta}{p}} |\Omega|^{\frac{\pt-\delta}{\pt}} \\
        & \leq \left( \left( \frac{1}{p} - \frac{1}{\pt} \right) \left( \frac{1}{\delta} - \frac{1}{\pt} \right)^{-1} \right) \|u\|_{\wop}^{p} + A \lambda^{\frac{p}{p-\delta}},
    \end{aligned}
\end{equation*}
where  
\[
A = \left[ \left( \frac{p}{\delta} \left( \frac{1}{p} - \frac{1}{\pt} \right) \left( \frac{1}{\delta} - \frac{1}{\pt} \right)^{-1} \right)^{\frac{-\delta}{p}} \|a\|_{L^{\infty}} S_{p}^{\frac{-\delta}{p}} |\Omega|^{\frac{\pt-\delta}{\pt}} \right]^{\frac{p}{p-\delta}}.
\]
This implies
    \begin{equation*}
        \jl(u) \geq - \left( \frac{1}{\delta} - \frac{1}{p_{*}} \right) A \lambda^{\frac{p}{p - \delta}},
    \end{equation*}
    which concludes the proof by setting $C_\delta$ as in \eqref{cdelta}. 
\end{proof}

\begin{lemma}\label{psrange}
    Let $\lambda \in (0, \lambda_0)$, and define $C_\delta$ as in \eqref{cdelta}. Furthermore $\{u_k\}\subset\nl$ be sequence such that  \( \jl(u_k) \to c \) and \( \jl^\prime(u_k) \to 0 \) as \( k \to \infty \). Then every such sequence has a convergent subsequence for values of $c$ satisfying
    \begin{equation*}
        -\infty < c < c_\infty := \frac{1}{N} \left( \frac{S_p}{\lambda} \right)^{\frac{N}{p}} - C_\delta \lambda^{\frac{p}{p-\delta}}.
    \end{equation*}
\end{lemma}
\begin{proof}
We have
    \begin{equation*}
        \begin{aligned}
            \frac{1}{p} \|u_k\|_{W_0^{1,p}}^p + \frac{1}{q} \|u_k\|_{\wsq}^q - \frac{\lambda}{\delta} \int_\Omega a(x) |u_k|^\delta - \frac{\lambda}{p_*} \int_\Omega |u_k|^{p_*} &= c + o_k(1), \\
            \|u_k\|_{W_0^{1,p}}^p + \|u_k\|_{\wsq}^q - \lambda \int_\Omega a(x) |u_k|^\delta - \lambda \int_\Omega |u_k|^{p_*} &= o_k(1).
        \end{aligned}
    \end{equation*}
Now $\{u_k\}$ is bounded in $W_0^{1,p}(\Omega)$, there exists $u \in W_0^{1,p}(\Omega)$ such that $u_k \rightharpoonup u$ weakly in $W_0^{1,p}(\Omega)$. Furthermore, $u$ is a critical point of $\jl$. We claim that $u_k \to u$ strongly in $W_0^{1,p}(\Omega)$. Since $u_k \to u$ strongly in $L^\gamma(\Omega)$ for all $1 \leq \gamma < p_*$, we obtain
    \begin{equation*}
        \int_\Omega a(x) |u_k|^\delta \to \int_\Omega a(x) |u|^\delta.
    \end{equation*}
Applying the Brezis-Lieb lemma, we obtain
    \begin{equation}\label{blatc}
        \begin{aligned}
            \frac{1}{p} \|u_k - u\|_{W_0^{1,p}}^p + \frac{1}{q} \|u_k - u\|_{\wsq}^q - \frac{\lambda}{p_*} \|u_k - u\|_{p_*}^{p_*} + \jl(u) \leq c + o_k(1).
        \end{aligned}
    \end{equation}
Additionally, we obtain the relation
    \begin{equation*}
        \|u_k - u\|_{W_0^{1,p}}^p + \|u_k - u\|_{\wsq}^q - \lambda \int_\Omega \left( |u_k|^{p_*} - |u|^{p_*} \right) = o_k(1).
    \end{equation*}
Defining $l = \lim_{k \to \infty} (\|u_k - u\|_{W_0^{1,p}}^p + \|u_k - u\|_{\wsq}^q)$, we deduce that $\lambda \int_\Omega (|u_k|^{p_*} - |u|^{p_*}) \to l$, which implies $\lambda \|u_k - u\|_{p_*}^{p_*} \to l$. If $l = 0$, then $u_k \to u$ strongly in $W_0^{1,p}(\Omega)$, completing the proof. Suppose instead that $l > 0$. Then
    \begin{equation*}
        l^{\frac{p}{p_*}} = \lambda \left( \lim_{k \to \infty} \int_\Omega |u_k - u|^{p_*} \right)^{\frac{p}{p_*}} \leq \lambda S_p^{-1} \lim_{k \to \infty} \|u_k - u\|_{W_0^{1,p}}^p \leq \lambda S_p^{-1} l.
    \end{equation*}
Consequently, we obtain $\left( \frac{S_p}{\lambda} \right)^{\frac{N}{p}} \leq l$. Using this and \eqref{blatc}, we establish
    \begin{equation*}
        \begin{aligned}
            c - \jl(u) &\geq \frac{1}{p} \|u_k - u\|_{W_0^{1,p}}^p + \frac{1}{q} \|u_k - u\|_{\wsq}^q - \frac{\lambda}{p_*} \|u_k - u\|_{p_*}^{p_*} + o_k(1) \\
            &\geq \left( \frac{1}{p} - \frac{1}{p_*} \right) l = \frac{l}{N}.
        \end{aligned}
    \end{equation*}
    This leads to
    \begin{equation*}
        c \geq \frac{l}{N} + \jl(u) \geq \frac{1}{N} \left( \frac{S_p}{\lambda} \right)^{\frac{N}{p}} - C_\delta \lambda^{\frac{p}{p-\delta}}.
    \end{equation*}
    Since this contradicts the assumption $c < c_\infty$, we conclude that $l = 0$, and thus $u_k \to u$ strongly in $W_0^{1,p}(\Omega)$. This completes the proof.
\end{proof}
\begin{theorem}\label{exsonesolcri}
There exists a constant $\Lambda_{0}>0$ such that for all $\lambda\in(0,\Lambda_{0})$, equation \eqref{generaleqn} admits a nontrivial nonnegative solution.
\end{theorem}
\begin{proof}
Define $\gamma_0>0$ such that for all $\lambda\in(0,\gamma_{0})$, the inequality
\begin{equation}\label{Lambda0}
    c_{\infty}= \frac{1}{N} \left(\frac{S_{p}}{\lambda}\right)^{\frac{N}{p}}-C_{\delta} \lambda^{\frac{p}{p-\delta}}>0 \text{ holds and set } \Lambda_{0}=\min\{\gamma_0,\lambda_{0}\}.
\end{equation}
 By Proposition \ref{criticallevel}, there exists a minimizing sequence $\{u_k\}$ for $\nl$ which is also  a $(PS)_{\theta_{\lambda}}$ sequence for $\jl$. Applying Lemma \ref{thetalambdaplus} and Lemma \ref{psrange}, we conclude that there exists $u_\lambda\in\wopo$ such that $u_k\to u_\lambda$ strongly in $\wopo$ for $\lambda\in(0,\Lambda_{0})$. Consequently, for such values of $\lambda$, the function $u_\lambda$ satisfies $\langle\jl^{\prime}(u_\lambda),u_\lambda\rangle=0,$  
and from \eqref{jleqninnl} it follows that
\begin{equation*}
    \jl(u_\lambda)\geq -\lambda\left(\frac{1}{\delta}-\frac{1}{\pt}\right)\int_{\Omega} a(x)|u|^{\delta}.
\end{equation*}
This leads to the estimate
\begin{equation*}
    \int_{\Omega} a(x)|u|^{\delta} \geq -\frac{\theta_{\lambda}}{\lambda}\left(\frac{1}{\delta}-\frac{1}{\pt}\right)^{-1}>0.
\end{equation*}
Thus, we conclude that $u_\lambda\not\equiv0$. This establishes that $u_\lambda\in\nl$ and satisfies $\jl(u_\lambda)=\theta_\lambda$. Next, we demonstrate that $u_\lambda\in\nl^{+}$. Suppose, for contradiction, that $u_\lambda\in\nl^{-}$. Then, from Lemma \ref{fiberingmapanalysis}, there exist $t_1<t_2=1$ such that $t_1 u_\lambda\in\nl^{+}$ and $t_2 u_\lambda\in\nl^{-}$. Since $\gamma_{u_{\lambda}}$ is increasing on $[t_1,t_2)$, it follows that
\begin{equation*}
    \theta_\lambda\leq \jl(t_1 u_\lambda)< \jl(t u_\lambda) \leq \jl( u_\lambda)=\theta_\lambda \quad \text{for } t\in(t_1,1),
\end{equation*}
which contradicts the assumption. Therefore, we conclude that $u_\lambda\in\nl^+$ and $\theta_\lambda=\jl(u_\lambda)=\theta_{\lambda}^{+}$.
Finally, using the same arguments as in the proof of Theorem \ref{subcrmul}, we establish that $u_\lambda$ is a nonnegative solution.
\end{proof}

We will establish the existence of a second solution below the first critical level using blowup analysis. To achieve this, we rely on asymptotic estimates of the minimizers of the Sobolev constant \( S_p \). The approach involves applying an appropriate truncation to the function:
\begin{equation*}
    U_{\varepsilon}(x)=\frac{K_{N,p}\varepsilon^{\frac{(N-p)}{p(p-1)}}}{(\varepsilon^{\frac{p}{p-1}}+|x|^{\frac{p}{p-1}})^{\frac{N-p}{p}}}, \quad \varepsilon>0, \quad K_{N,p}=\left[N\left(\frac{N-p}{p-1}\right)^{p-1}\right]^{\frac{N-p}{p^2}}.
\end{equation*}
Clearly, \( U_{\varepsilon} \in W_0^{1,p}(\mathbb{R}^{N}) \) and attains the best Sobolev constant \( S_p \). Let us fix \( r>0 \) such that \( \overline{B_{4r}(0)} \subset \Omega \) and introduce a radial cutoff function \( \phi_{r} \in C^{\infty}(\mathbb{R}^{N},[0,1]) \) satisfying:
\begin{equation*}
    \phi_{r}=
    \begin{cases}
        1, & x \in B_{r}, \\
        0, & x \in B_{2r}^{c}.
    \end{cases}
\end{equation*}
For any \( \varepsilon>0 \), we define \( u_\varepsilon = \phi_{r} U_{\varepsilon} \) and \( v_\varepsilon = \frac{u_\varepsilon}{\|u_\varepsilon\|_{p_*}} \), both of which are radial and belong to \( W_0^{1,p}(\mathbb{R}^N) \). From \cite{GA87} and \cite[(3.6)-(3.9)]{GA91}, for sufficiently small \( \varepsilon \), we obtain:
\begin{align}
    \|v_\varepsilon\|_{W_0^{1,p}}^{p} &= S_{p} + O(\varepsilon^{\frac{N-p}{p-1}}), \label{sobnorm} \\
    C_1 \varepsilon^{N-\frac{t(N-p)}{p}} &\leq \|v_\varepsilon\|_{L^t}^{t} \leq C_2 \varepsilon^{N-\frac{t(N-p)}{p}}, \quad \text{for } t>p_*(1-\frac{1}{p}), \label{tlargenorm} \\
    C_1 \varepsilon^{(N-p)t/p^2}|\log \varepsilon| &\leq \|v_\varepsilon\|_{L^t}^{t} \leq C_2 \varepsilon^{(N-p)t/p^2}|\log \varepsilon|, \quad \text{for } t=p_*(1-\frac{1}{p}), \label{tmiddlenorm} \\
    C_1 \varepsilon^{\frac{(N-p)t}{p(p-1)}} &\leq \|v_\varepsilon\|_{L^t}^{t} \leq C_2 \varepsilon^{\frac{(N-p)t}{p(p-1)}}, \quad \text{for } t<p_*(1-\frac{1}{p}). \label{tsmallnorm}
\end{align}
We now proceed to estimate \( \| v_\varepsilon \|_{\wsq} \).
\begin{lemma}
    \label{wsqnorm}
    Define \( \mpqcr = \min \left\{ \frac{q(N-p)}{p(p-1)}, q(1-s) + N \left( 1 - \frac{q}{p} \right) \right\}. \) Then, for sufficiently small \( \varepsilon \), we have
    \begin{equation*}
        \| u_\varepsilon \|_{\wsq}^{q} = O(\varepsilon^\mpqcr).
    \end{equation*}
\end{lemma}
\begin{proof}
 We partition \( \mathbb{R}^{2N} \) into four subdomains \( D_{i}, i=1,2,3,4 \), where \( \bigcup_{i=1}^{4} D_{i} = \mathbb{R}^{2N} \):
    \begin{equation*}
        \begin{aligned}
            D_{1} &= B_{r} \times B_{r}, \\
            D_{2} &= \{(x,y) \in B_{r} \times B_{r}^{c} : |x-y| > r/2 \}, \\
            D_{3} &= \{(x,y) \in B_{r} \times B_{r}^{c} : |x-y| \leq r/2 \}, \\
            D_{4} &= B_{r}^{c} \times B_{r}^{c}.
        \end{aligned}
    \end{equation*}
From \cite[(5.6) and (5.7)]{DFV24}, for sufficiently small \( \varepsilon \), we obtain:
    \begin{align*}
        | u_\varepsilon (x) - u_\varepsilon (y) | &\leq C \varepsilon^{\frac{(N-p)}{p(p-1)}} |x-y|, \quad \text{if } x \in \mathbb{R}^N, y \in B_{r}^{c}, |x-y| \leq r/2, \\
        | u_\varepsilon (x) - u_\varepsilon (y) | &\leq C \varepsilon^{\frac{(N-p)}{p(p-1)}} \min\{1,|x-y|\}, \quad \text{if } x, y \in B_{r}^{c}.
    \end{align*}
 Using these estimates we infer that 
    \begin{equation*}
        \begin{aligned}
           & \underset{D_{4}}{\int\int} \frac{|\ueps(x)-\ueps(y)|^{q}}{|x-y|^{N+sq}} \leq C \varepsilon^\frac{q(N-p)}{p(p-1)} \underset{B_{2r} \times \RR^N}{\int\int} \frac{1}{|x-y|^{N+sq}} \leq O(\varepsilon^\frac{q(N-p)}{p(p-1)}),
           \\
           & \underset{D_{3}}{\int\int} \frac{|\ueps(x)-\ueps(y)|^{q}}{|x-y|^{N+sq}} \leq C \varepsilon^\frac{q(N-p)}{p(p-1)}  \underset{D_3}{\int \int} \frac{1}{|x-y|^{N+sq}}
           \\
           &\hspace{5cm}\leq C \varepsilon^\frac{q(N-p)}{p(p-1)} \int_{B_{r}} dx \int_{B_{r/2}} \frac{1}{\zeta^{N+sq-q}} d\zeta = C \varepsilon^\frac{q(N-p)}{p(p-1)} =O(\varepsilon^\frac{q(N-p)}{p(p-1)}).
        \end{aligned}
    \end{equation*}
   We have that $\ueps=U_{\varepsilon}$ in $B_{r}$ and 
    \begin{equation*}
        \begin{aligned}
            |\ueps(x)-\ueps(y)|^{q} \leq 2^{q-1}(|U_\varepsilon(x)-U_\varepsilon(y)|^{q}+|U_\varepsilon(y)-\ueps(y)|^{q}) \text{ in } D_2.
        \end{aligned}
    \end{equation*}
    Thus we get that
    \begin{equation*}
        \begin{aligned}
            \underset{D_{2}}{\int\int} \frac{|U_\varepsilon(y)-\ueps(y)|^{q}}{|x-y|^{N+sq}} \thinspace dxdy &\leq \underset{D_{2}}{\int\int} \frac{2^{q-1}(|U_\varepsilon(y)|^{q}+|\ueps(y)|^{q}}{|x-y|^{N+sq}} \thinspace dxdy \leq \underset{D_{2}}{\int\int} \frac{2^{q}|U_\varepsilon(y)|^{q}}{|x-y|^{N+sq}} \thinspace dxdy 
            \\
            &\leq  C \varepsilon^\frac{q(N-p)}{p(p-1)} \int_{B_{r}} dx \int_{B_{r/2}^{c}} \frac{1}{|\zeta|^{N+sq}} d\zeta = C \varepsilon^\frac{q(N-p)}{p(p-1)}.
        \end{aligned}
    \end{equation*}
   Observe that if $x=\varepsilon\xi$
   \begin{equation*}
       \begin{aligned}
                   U_{\varepsilon}(x)=\frac{K_{N,p}\varepsilon^{\frac{(N-p)}{p(p-1)}}}{(\varepsilon^{\frac{p}{p-1}}+|x|^{\frac{p}{p-1}})^{\frac{N-p}{p}}} = \varepsilon^{\frac{N-p}{p(p-1)}-\frac{N-p}{p-1}}U_{1}(\xi)=\varepsilon^{-\frac{N-p}{p}} U_{1}(\xi).
       \end{aligned}
   \end{equation*}
   Finally, by a double change of variables $x=\varepsilon\xi$ and $y=\varepsilon\varsigma$,
    \begin{equation*}
        \begin{aligned}
              \underset{\RR^{2N}}{\int\int} \frac{|U_\varepsilon(x)-U_\varepsilon(y)|^{q}}{|x-y|^{N+sq}} \thinspace dxdy
              =\varepsilon^{-\frac{q}{p}(N-p)+N-sq}  \underset{\RR^{2N}}{\int\int} \frac{|U_1(\xi)-U_1(\varsigma)|^{q}}{|\xi-\varsigma|^{N+sq}} \leq C \varepsilon^{N(1-\frac{q}{p})+q(1-s)}.
        \end{aligned}
    \end{equation*}
    Hence for any $V\subset\RR^{2N},$ in particular for $V=D_1,D_2$ we have  
    \begin{equation*}
        \begin{aligned}
            \underset{V}{\int\int} \frac{|U_\varepsilon(x)-U_\varepsilon(y)|^{q}}{|x-y|^{N+sq}} \thinspace dxdy \leq C \varepsilon^{N(1-\frac{q}{p})+q(1-s)}.
        \end{aligned}
    \end{equation*}
    Now combining all the estimate we conclude that $\|\ueps\|_{\wsq}^{q}\leq\varepsilon^\mpqcr.$
\end{proof}
Using the estimate \( \| u_\varepsilon \|_{p_*} = O(1) \), we derive the following corollary:
\begin{corollary}
    Let \( \mpqcr\) be as defined in Lemma \ref{wsqnorm}. Then, for sufficiently small \( \varepsilon \),
    \begin{equation*}
        \| v_\varepsilon \|_{\wsq}^{q} = O(\varepsilon^\mpqcr).
    \end{equation*}
\end{corollary}
We aim to establish the existence of a second solution lying below the first critical level.
\begin{lemma}\label{talentili}
    There exists a constant $\Lambda_{0} > 0$ such that for every $\lambda \in (0, \Lambda_{0})$, there exists a function $u \geq 0$ in $\wopo$ satisfying
    \begin{equation*}
        \sup_{t \geq 0} \jl(tu) < c_{\infty}.
    \end{equation*}
In particular, this implies that $\theta_{\lambda}^{-} < c_{\infty}$.
\end{lemma}
\begin{proof}
   Without loss of generality, we assume that $\inf_{B_{r_0}} a(x) = m_{a} > 0.$ Let $\Lambda_0$ be as defined in \eqref{Lambda0}, ensuring that $c_{\infty} > 0$ for all $\lambda \in (0, \Lambda_0)$. Consider the function $\veps$ and get 
    \begin{equation*}
        \jl(t\veps) \leq \frac{t^{p}}{p} \|\veps\|_{\wop}^{p} + \frac{t^{q}}{q} \|\veps\|_{\wsq}^{q} \leq C (t^{p} + t^{q}).
    \end{equation*}
    Consequently, there exists some $t_0 \in (0,1)$ such that
    \begin{equation*}
        \sup_{0 \leq t \leq t_0} \jl(t\veps) < c_{\infty}.
    \end{equation*}
    Next, define the function
    \begin{equation*}
        h(t) = \frac{t^p}{p}\|\veps\|_{\wop}^{p} + \frac{t^q}{q} \|\veps\|_{\wsq}^{q} - \lambda\frac{t^{\pt}}{\pt}.
    \end{equation*}
    Noting that $h(0) = 0$, that $h(t) > 0$ for small $t$, and that $h(t) < 0$ for sufficiently large $t$, there exists $t_\varepsilon > 0$ such that
    \begin{equation*}
        \max_{t > 0} h(t) = h(t_\varepsilon),
    \end{equation*}
    where $t_\varepsilon$ is determined by solving $h'(t_\varepsilon) = 0$, leading to
   \begin{equation*}
        \begin{aligned}
            t_\varepsilon&= \frac{1}{\lambda^{\frac{1}{\pt-q}}}(t_\varepsilon^{p-q}\|\veps\|_{\wop}^{p}+O(\varepsilon^\mpqcr) )^{\frac{1}{\pt-q}}.
            \end{aligned}
\end{equation*}
 We choose $\varepsilon^{\frac{1}{\beta}}=\lambda$ for $\beta>\frac{1}{\mpqcr}$ to infer that 
\begin{equation}\label{tepseqn}
        \begin{aligned}
          t_\epsilon \lesssim  \; t_\varepsilon^{\frac{p-q}{\pt-q}} \left(\frac{\|\veps\|_{\wop}^p}{\lambda}\right)^\frac{1}{\pt-q}+ \varepsilon^{\frac{\mpqcr}{\pt-q}-\frac{1}{\beta(\pt-q)}}.
        \end{aligned}
    \end{equation}
   We claim that : 
    \begin{equation*}
        t_\varepsilon \lesssim \left(\frac{\|\veps\|_{\wop}^p}{\lambda}\right)^{\frac{1}{\pt - p}}.
    \end{equation*}
   We first note that if $t_\varepsilon$ is uniformly bounded then claim is true as $\frac{\|\veps\|_{\wop}^p}{\lambda}\rightarrow \infty$ as $\varepsilon\rightarrow 0.$ If $t_\varepsilon$ is not bounded, then from \eqref{tepseqn} , we can write
   \begin{equation*}
       t_\varepsilon \leq 2 t_{\varepsilon}^{\frac{p-q}{\pt-q}} \left(\frac{\|\veps\|_{\wop}^p}{\lambda}\right)^\frac{1}{\pt-q}
   \end{equation*}
   and this proves the claim.\\ 
   Next, to find an upper bound for $\jl(t\veps)$ we use the estimate 
    \begin{equation*}
        \int_{\Omega} a(x) |\veps|^{\delta} \geq m_{a} \int_{B_{r_0}} |\veps|^{\delta},
    \end{equation*}
    and obtain 
    \begin{equation*}
    \begin{aligned}
            \sup_{t \geq t_0} \jl(t\veps) &\leq \sup_{t > 0} h(t) - \lambda \frac{t_0^\delta m_a}{\delta} \int_{B_{r_0}} |\veps|^{\delta} \\
            &= \frac{t_\varepsilon^p}{p}\|\veps\|_{\wop}^{p}+ \frac{t_\varepsilon^q}{q} \|\veps\|_{\wsq}^{q}-\lambda\frac{t_\varepsilon^{\pt}}{\pt} -\lambda\frac{t_0^\delta m_a}{\delta}\int_{B_{r_0}} |\veps|^{\delta}\\
        & \leq \sup_{t\geq 0}\left( \frac{t^p}{p}\|\veps\|_{\wop}^{p} -\lambda\frac{t^{\pt}}{\pt}\right)+ \frac{\left( \|\veps\|_{\wop}^p\right)^\frac{q}{\pt-p}}{\lambda^\frac{q}{\pt-p} q} O(\varepsilon^\mpqcr)-\lambda\frac{t_0^\delta m_a}{\delta}\int_{B_{r_0}} |\veps|^{\delta}.
    \end{aligned}
    \end{equation*}
    Defining $g(t) = \frac{t^p}{p} \|\veps\|_{\wop}^{p} - \lambda \frac{t^{\pt}}{\pt}$, we observe that $g$ attains its maximum at $\widetilde{t} = \left(\frac{\|\veps\|_{\wop}^{p}}{\lambda}\right)^{\frac{1}{\pt - p}}$. Consequently,
    \begin{equation*}
        \sup_{t \geq 0} g(t) = g(\widetilde{t}) = \frac{1}{N} \left(\frac{\|\veps\|_{\wop}^N}{\lambda^{(N - p)/p}}\right) \leq \frac{1}{\lambda^{(N - p)/p} N} \left(S_{p}^{N/p} + O(\varepsilon^{(N - p)/(p - 1)})\right).
    \end{equation*}
    This leads to the final estimate
    \begin{equation*}
        \sup_{t \geq t_{0}} \jl(t\veps) \leq \frac{1}{N}\left(\frac{S_p}{\lambda}\right)^{N/p}+\frac{S_p^{N/p}}{N\lambda^{N/p}}(\lambda-1) + \frac{1}{\lambda^{\frac{N-p}{p}}}O(\varepsilon^\mpqcr)-\lambda\int_{B_{r_0}} |\veps|^{\delta}.
    \end{equation*}
    Set $\varepsilon = \lambda^\beta$ for $\beta \geq \max\left\{\frac{1}{\mpqcr} \left(\frac{N}{p} - 1\right),\frac{1}{\mpqcr}\right\}$, and obtain
    \begin{equation*}
        \sup_{t \geq t_{0}} \jl(t\veps) \leq \frac{1}{N} \left(\frac{S_p}{\lambda}\right)^{N/p} - C_{\delta} \lambda^{p/(p - \delta)}
    \end{equation*}
    for sufficiently small $\lambda.$ By applying the Lemma \ref{fiberingmapanalysis}, we deduce the existence of $\hat{t} > 0$ such that $\hat{t} \veps \in \nl^{-}$, thus concluding that $\theta_\lambda^- < c_{\infty}$.
\end{proof}
Now we prove the existence of second solution for the critical nonlinearity.
\\[1mm]
\textit{\underline{\textbf{Proof of Theorem \ref{multsolcri}} }}
\\[3mm]
    The results of  Proposition \ref{criticallevel} remain valid even when $\nl$ is replaced by $\nl^-$. Consequently, we obtain a minimizing sequence $\{u_k\} \subset \nl^-$ satisfying
    \begin{equation*}
        \jl(u_k) = \theta_\lambda^- + o_k(1), \quad \text{and} \quad \jl'(u_k) = o_k(1),
    \end{equation*}
    implying that $\{u_k\}$ forms a $(PS)_{\theta_\lambda^-}$ sequence for $\jl$.
    By applying Lemmas \ref{psrange} and \ref{talentili}, there exists a function $v_\lambda \in \wopo$ such that $u_k \to v_\lambda$ in $\wopo$. Theorem \ref{cdeltathm} further guarantees that $(\jl'(v_\lambda), v_\lambda) = 0$.
    Exploiting the strong convergence $u_k \to v_\lambda$ and observing that $\nl^0 = \varnothing$, we deduce that $v_\lambda \in \nl^-$ and $\theta_\lambda^- = \jl(v_\lambda)$. Finally, following an argument similar to the proof of Theorem \ref{subcrmul}, we conclude that $v_\lambda$ is non-negative. \hfill\qed
\section{Brezis Nirenberg Type Problem}\label{brezisnirenberg}
In this section, we examine the scenario where \( q \leq p \), \( r = \pt \), and \( b(x) = \frac{1}{\lambda} \) in equation \eqref{generaleqn}. Under these conditions, equation \eqref{generaleqn} reduces to  
\begin{equation}\label{bneqn}
    \begin{aligned}
        \fp u + \fqs u &= \lambda a(x)|u|^{\delta-2}u + |u|^{\pt-2}u \quad \text{in } \Omega, \\
        u &= 0 \quad \text{in } \Omega^{c}.
    \end{aligned}
\end{equation}
Our focus is on establishing the existence of multiple nonnegative solutions for this class of problems. The associated energy functional \( \Jl : \mathcal{W}_0(\Omega) \to \mathbb{R} \) is given by
\begin{equation*}
    \Jl(u) := \frac{1}{p} \|u\|_{\wop}^{p} + \frac{1}{q} \|u\|_{\wsq}^{q} -  \int_{\Omega} \left( \lambda\frac{a(x)}{\delta} |u|^{\delta} + \frac{b(x)}{\pt} |u|^{\pt} \right) dx.
\end{equation*}  
The Nehari manifold corresponding to \( \Jl \) is defined as  
\begin{equation*}
    \begin{aligned}
        \Nl &:= \{u \in \mathcal{W}_0(\Omega) \setminus \{0\} : \langle \Jl^{\prime}(u), u \rangle = 0\}.
    \end{aligned}
\end{equation*}
We define the fibering map associated with \( \Jl \) as \( \Gamma_{u}:\mathbb{R}_{+}\to\mathbb{R} \) by  $\Gamma_u(t) = \Jl(tu).$ Depending on the behavior of $\Gamma_u,$ we set $\Nl^0,\Theta_\lambda, \Nl^{\pm}$ and $\Theta_\lambda^{\pm}$ similar to section \ref{nehari}. Define the functional \( \mathcal{E}_{\lambda}: \Nl \to \mathbb{R} \) as  
\begin{equation*}
    \mathcal{E}_{\lambda}(u) := \frac{\pt-p}{\pt-\delta} \|u\|_{\wop}^{p} + \frac{\pt-q}{\pt-\delta} \|u\|_{\wsq}^{q} - \lambda \int_{\Omega} a(x) |u|^{\delta}.
\end{equation*}
Similar to subsection \ref{fibermapsec}, we also want to give complete characterization of the geometry of the fibering maps associated with problem \eqref{bneqn}. To this end, we introduce the auxiliary $C^{1}$ function  $M_{u}:\RR_{+}\rightarrow\RR$ which is defined for a fixed $u\in\wo\setminus\{0\}$ as
\begin{equation*}
    \begin{aligned}
        M_{u}(t)= t^{(p-\delta)} \|u\|_{\wop}^{p}+ t^{(q-\delta)}\|u\|_{\wsq}^{q}- t^{(\pt-\delta)} \int_{\Omega}  |u|^{\pt}\text{ for } t\geq 0. 
    \end{aligned}
\end{equation*} 
It follows that $tu\in\Nl$ if and only if $t$ satisfies $M_{u}(t)=\lambda  \int_{\Omega}  a(x)|u|^{\delta}.$ If $\int_{\Omega}a(x)|u|^{\delta}>0,$  we see $ M_{u}(t)\rightarrow -\infty$ as $t\rightarrow\infty,$ $ M_{u}(t)>0$ for $t$ small enough and $ M_{u}^{\prime}(t)<0$ for $t$ large enough. Following a similar argument to subsection \ref{fibermapsec}, there exists unique $t_{max}>0$ such that $ M_{u}^{\prime}(t_{max})=0.$ Furthermore there exist \( t_{1} < t_{\max} \) and \( t_{2} > t_{\max} \) such that \( t_{1}u \in \Nl^{+} \) and \( t_{2}u \in \Nl^{-} \). Additionally we also get $\Jl(t_{1}u) = \min_{t \in [0, t_{2}]} \Jl(tu)$ and $\Jl(t_2 u) = \max_{t \geq t_{\max} } \Jl (t u).$  Now suppose \( \int_{\Omega} a(x) |u|^{\delta} < 0 \).  
We observe that \( M_{u}(t) \to -\infty \) as \( t \to \infty \), while \( M_{u}(t) > 0 \) for sufficiently small \( t \), and \( M_{u}^{\prime}(t) < 0 \) for large \( t \).  
By an argument similar to subsection \ref{fibermapsec}, there exists a unique \( t_0 > 0 \) such that \( M_{u} \) is increasing on \( (0,t_0) \) and decreasing on \( (t_0,\infty) \), with \( M_{u}^{\prime}(t_0) = 0 \).  
Since \( M_{u}(t_0) > 0 \) and \( \lambda \int_{\Omega} a(x) |u|^{\delta} < 0 \), there exists a unique \( t_1 > 0 \) such that  
\[
M_{u}(t_1) = \lambda \int_{\Omega} a(x) |u|^{\delta}, \quad M_{u}^{\prime}(t_1) < 0.
\]
This implies that \( t_1 u \in \Nl^{-} \), meaning \( t_1 u \) is a local maximum. Now similar to Lemma \ref{thetalambdaplus} we get the following lemma 
\begin{lemma}\label{thetalambdaplusbn}
   There exists constant \( C_2 > 0 \) such that 
   \begin{equation*}
  \Theta_{\lambda}^{+}\leq   \begin{aligned}
& -\frac{(p-\delta)(\pt-p)}{p\delta \pt} C_2 <0.
 \end{aligned}
\end{equation*}    
\end{lemma} 
Similar to Lemma \ref{nlambdazero}, there exists  \( \bar{\lambda}_{0} > 0 \) such that \( \Nl^{0} = \varnothing \) for all \( \lambda \in (0, \bar{\lambda}_{0}) \). For \( z \in \Nl \), define the function \( \mathcal{H}_z: \mathbb{R} \times \wo(\Omega) \to \mathbb{R} \) by 
\begin{equation*}
    \begin{aligned}
        \mathcal{H}_{z}(t,w) &= \langle \Jl^{\prime}(t(z-w)), t(z-w) \rangle.
    \end{aligned}
\end{equation*}
Using arguments similar to Lemma \ref{projlem} we get the following lemma: 
\begin{lemma}\label{projlembn}
   Let \( \lambda \in (0, \bar{\lambda}_{0}) \) and \( z \in \Nl \). Then there exist \( \varepsilon > 0 \) and a differentiable function \( \Pi: B(0, \varepsilon) \subset \wo(\Omega) \to \mathbb{R}_{+} \) such that \( \Pi(0) = 1 \), \( \Pi(w)(z-w) \in \Nl \), and 
   \begin{equation}
       \begin{aligned}
           \langle \Pi^{\prime}(0),w \rangle = 
           \frac{p \mathfrak{A}_{p}(z,w) + q A_{q}(z,w,\mathbb{R}^{2N}) -  \int_{\Omega} \left( \lambda\delta a(x) |z|^{\delta-1}w + \pt  |z|^{\pt-1}w \right)}
           {(p-\delta) \|z\|_{W_0^{1,p}}^{p} + (q-\delta) \|z\|_{W_0^{s,q}}^{q} -  (\pt-\delta) \int_{\Omega}  |z|^{\pt} }.
       \end{aligned}
   \end{equation}
\end{lemma}
Our next objective is to establish the existence of a sequence that satisfies the $(PS)$ condition. 
\begin{proposition}\label{criticallevelbn}
   Let  $\lambda \in (0,\bar{\lambda}_{0})$. Then, there exists a sequence $\{u_k\} \subset \Nl$ such that
   \begin{equation*}
       \Jl(u_k) = \Theta_\lambda + o_k(1), \quad \Jl^{\prime}(u_k) = o_k(1).
   \end{equation*}
\end{proposition}
\begin{proof}
    Since $\Jl$ is coercive and bounded below in $\Nl$, the Ekeland variational principle guarantees the existence of a minimizing sequence $\{u_k\} \subset \Nl$ satisfying
   \begin{equation}\label{ekelandseqbn}
       \Jl(u_k) < \Theta_\lambda + \frac{1}{k}, \quad \Jl(u_k) < \Jl(v) + \frac{1}{k} \|v - u_k\|_{\wo}, \quad \forall v \in \Nl.
   \end{equation}
   Since $u_k \in \Nl$, we have
   \begin{equation*}
       \Jl(u_k) = \left(\frac{1}{p} - \frac{1}{\pt}\right) \|u_k\|_{\wop}^{p} + \left(\frac{1}{q} - \frac{1}{\pt}\right) \|u_k\|_{\wsq}^{q} - \lambda \left(\frac{1}{\delta} - \frac{1}{\pt}\right) \int_\Omega a(x) |u_k|^{\delta}.
   \end{equation*}
   By \eqref{ekelandseqbn} and Lemma \ref{thetalambdaplusbn},
   \begin{equation*}
    \begin{aligned}
        \Jl(u_k) <\Theta_\lambda+ \frac{1}{k}<\Theta_\lambda^{+}<0.
    \end{aligned}
\end{equation*} 
and it follows that $u_k \not\equiv 0$ for sufficiently large $k$. Moreover, applying H\"older's inequality, we obtain
\begin{equation*}
    \begin{aligned}
        \left(\frac{(-\Theta_\lambda^{+})\delta r S_{p}^{\frac{\delta}{p}}}{(r-\delta)\anorm\lambda}\right)^{\frac{1}{\delta}} \leq  \|u_k\|_{\wop}\leq \left(\frac{\lambda(r-\delta)p\anorm}{(r-p)\delta S_{p}^{\frac{\delta}{p}}}\right)^{\frac{1}{p-\delta}}
    \end{aligned}
\end{equation*}
Next, we aim to prove that $\|\Jl^{\prime}(u_k)\|\rightarrow 0$ as $k\to\infty$. Using
Lemma \ref{projlembn}, for each $u_k$, there exist $\varepsilon_k>0$ and differentiable functions $\Pi_k: B(0, \varepsilon_k) \subset \wo \to \mathbb{R}_+$ with $\Pi_k(0) = 1$ and $\Pi_k(v)(u_k - v) \in \nl$ for all $v \in B(0, \varepsilon_k)$. Fix $k\in\mathbb{N}$ such that $u_{k}\not\equiv 0$ and $0<\rho<\varepsilon_k.$ Set $v_\rho = \frac{\rho u_k}{\|u_k\|_{\wop}}$ and $h_\rho = \Pi_k(v_\rho)(u_k - v_\rho).$ We obtain the estimate  
\begin{equation*}  
    \langle \Jl^{\prime}(u_k), \frac{u_k}{\|u_k\|_{\wop}} \rangle \leq \frac{C}{k} (1 + \|\xi_k^{\prime}(0)\|_{\wop}),  
\end{equation*}  
following a similar argument as in Lemma \ref{projlembn}. Furthermore, from \eqref{ekelandseqbn}, we deduce that $\|\xi_k^{\prime}(0)\|_{\wo}$ is bounded. This completes the proof.
\end{proof}

\begin{theorem}\label{cdeltathmbn}
    Let $\lambda \in (0, \bar{\lambda}_0)$, and suppose that $\{u_k\} \subset \Nl$ is a $(PS)_c$ sequence for $\Jl$, with $u_k$ converging weakly to $u$ in $\wop$. Then, $\Jl^\prime(u)=0$ and there exists a positive constant $C_{\delta}$, depending on $p, N, S_p, |\Omega|$, and $\delta$, such that
    \begin{equation}
        \Jl(u) \geq -C_{\delta} \lambda^{\frac{p}{p - \delta}},
    \end{equation}
    where
    \begin{equation}\label{cdeltabn}
        C_{\delta} = \left( \frac{1}{\delta} - \frac{1}{p_{*}} \right) \left[ \left( \frac{p}{\delta} \left( \frac{1}{p} - \frac{1}{p_{*}} \right) \left( \frac{1}{\delta} - \frac{1}{p_{*}} \right)^{-1} \right)^{-\frac{\delta}{p}} \|a\|_{L^{\infty}} S_p^{-\frac{\delta}{p}} |\Omega|^{\frac{p_{*} - \delta}{p_{*}}} \right]^{\frac{p}{p - \delta}}.
    \end{equation}
\end{theorem}
\begin{sketch}

We argue similar to Theorem \ref{cdeltathm} and that leads us to
    \begin{equation*}
        \Jl(u) \geq \left( \frac{1}{p} - \frac{1}{p_{*}} \right) \|u\|_{\wop}^{p} - \lambda \left( \frac{1}{\delta} - \frac{1}{p_{*}} \right) \int_{\Omega} a(x) |u(x)|^{\delta} \,dx.
    \end{equation*}
    Applying H\"older's inequality, Sobolev embeddings, and Young’s inequality, we derive
  \begin{equation*}
        \begin{aligned}
            \lambda \int_{\Omega} a |u|^\delta &= (\frac{p}{\delta}(\frac{1}{p}-\frac{1}{\pt})(\frac{1}{\delta}-\frac{1}{\pt})^{-1})^{\frac{\delta}{p}} \|u\|_{\wop}^{\delta} \times \lambda (\frac{p}{\delta}(\frac{1}{p}-\frac{1}{\pt})(\frac{1}{\delta}-\frac{1}{\pt})^{-1})^{\frac{-\delta}{p}} \|a\|_{L^{\infty}} S_{p}^{\frac{-\delta}{p}} |\Omega|^{\frac{\pt-\delta}{\pt}}
            \\
            & \leq ((\frac{1}{p}-\frac{1}{\pt})(\frac{1}{\delta}-\frac{1}{\pt})^{-1}) \|u\|_{\wop}^{p} + A \lambda^{\frac{p}{p-\delta}}
        \end{aligned}
    \end{equation*}
    where $A$ is given in \eqref{cdeltabn}. This implies
    \begin{equation*}
        \Jl(u) \geq - \left( \frac{1}{\delta} - \frac{1}{p_{*}} \right) A \lambda^{\frac{p}{p - \delta}},
    \end{equation*}
    which concludes the proof.
\end{sketch}
\begin{lemma} \label{psrangebn}
    Let $\lambda \in (0, \bar\lambda_0)$, and define $C_\delta$ as in \eqref{cdeltabn}. Furthermore $\{u_k\}\subset\Nl$ be sequence such that  \( \Jl(u_k) \to c \) and \( \Jl^\prime(u_k) \to 0 \) as \( k \to \infty \). Then every such sequence
        \begin{equation*}
        -\infty < c < C_\infty := \frac{S_p^{\frac{N}{p}}}{N}  - C_\delta \lambda^{\frac{p}{p-\delta}}.
    \end{equation*}
    \end{lemma}
\begin{proof}
Consider a $(PS)_c$ sequence $\{u_k\}$ for $\Jl$ in $W_0^{1,p}(\Omega)$. Then we have
    \begin{equation*}
        \begin{aligned}
            \frac{1}{p} \|u_k\|_{W_0^{1,p}}^p + \frac{1}{q} \|u_k\|_{\wsq}^q - \frac{\lambda}{\delta} \int_\Omega a(x) |u_k|^\delta - \frac{1}{p_*} \int_\Omega |u_k|^{p_*} &= c + o_k(1), \\
            \|u_k\|_{W_0^{1,p}}^p + \|u_k\|_{\wsq}^q - \lambda \int_\Omega a(x) |u_k|^\delta -  \int_\Omega |u_k|^{p_*} &= o_k(1).
        \end{aligned}
    \end{equation*}
Now $\{u_k\}$ is bounded in $W_0^{1,p}(\Omega)$, there exists $u \in W_0^{1,p}(\Omega)$ such that $u_k \rightharpoonup u$ weakly in $W_0^{1,p}(\Omega)$. Furthermore, $u$ is a critical point of $\Jl$. We claim that $u_k \to u$ strongly in $W_0^{1,p}(\Omega)$. Since $u_k \to u$ strongly in $L^\gamma(\Omega)$ for all $1 \leq \gamma < p_*$, we obtain
    \begin{equation*}
        \int_\Omega a(x) |u_k|^\delta \to \int_\Omega a(x) |u|^\delta.
    \end{equation*}
Applying the Brezis-Lieb lemma, we obtain
    \begin{equation*}
        \begin{aligned}
            \frac{1}{p} \|u_k - u\|_{W_0^{1,p}}^p + \frac{1}{q} \|u_k - u\|_{\wsq}^q - \frac{1}{p_*} \|u_k - u\|_{p_*}^{p_*} + \Jl(u) \leq c + o_k(1).
        \end{aligned}
    \end{equation*}
Additionally, we obtain the relation
    \begin{equation*}
        \|u_k - u\|_{W_0^{1,p}}^p + \|u_k - u\|_{\wsq}^q -  \int_\Omega \left( |u_k|^{p_*} - |u|^{p_*} \right) = o_k(1).
    \end{equation*}
Defining $l = \lim_{k \to \infty} (\|u_k - u\|_{W_0^{1,p}}^p + \|u_k - u\|_{\wsq}^q)$, we deduce that $ \int_\Omega (|u_k|^{p_*} - |u|^{p_*}) \to l$, which implies $ \|u_k - u\|_{p_*}^{p_*} \to l$. If $l = 0$, then $u_k \to u$ strongly in $W_0^{1,p}(\Omega)$, completing the proof. Suppose instead that $l > 0$. Then
    \begin{equation*}
        l^{\frac{p}{p_*}} =  \left( \lim_{k \to \infty} \int_\Omega |u_k - u|^{p_*} \right)^{\frac{p}{p_*}} \leq  S_p^{-1} \lim_{k \to \infty} \|u_k - u\|_{W_0^{1,p}}^p \leq  S_p^{-1} l.
    \end{equation*}
Consequently, we obtain $ S_p^{\frac{N}{p}} \leq l$. Using this, we establish
    \begin{equation*}
        \begin{aligned}
            c - \Jl(u) &\geq \frac{1}{p} \|u_k - u\|_{W_0^{1,p}}^p + \frac{1}{q} \|u_k - u\|_{\wsq}^q - \frac{1}{p_*} \|u_k - u\|_{p_*}^{p_*} + o_k(1) \\
            &\geq \left( \frac{1}{p} - \frac{1}{p_*} \right) l = \frac{l}{N}.
        \end{aligned}
    \end{equation*}
    This leads to
    \begin{equation*}
        c \geq \frac{l}{N} + \Jl(u) \geq \frac{S_p^{\frac{N}{p}}}{N} - C_\delta \lambda^{\frac{p}{p-\delta}}.
    \end{equation*}
    Since this contradicts the assumption $c < C_\infty$, we conclude that $l = 0$, and thus $u_k \to u$ strongly in $W_0^{1,p}(\Omega)$. This completes the proof.
\end{proof}
Define $\mu_0>0$ such that for all $\lambda\in(0,\mu_{0})$, the inequality
\begin{equation}\label{Lambda0bn}
    C_{\infty}\geq \frac{S_{p}^{\frac{N}{p}}}{N} -C_{\delta} \lambda^{\frac{p}{p-\delta}}>0 \text{ holds and set } \bar\Lambda_{0}=\min\{\mu_0,\bar\lambda_{0}\}.
\end{equation}

\begin{lemma}\label{talentilibn}
Assume that either
\begin{equation*}
    \max\left\{\frac{Np}{\mpqcr+N-p},\pt(1-1/p)\right\}<\delta<q
\end{equation*} or 
\begin{equation*}
    \delta<\pt(1-1/p) \text{ and } s<1-\frac{1}{q}\left(\frac{N-p}{p-1}-N\left(1-\frac{q}{p}\right)\right).
\end{equation*}      Then there exists a constant $\bar\Lambda_{00} > 0$ such that for every $\lambda \in (0, \bar\Lambda_{00})$, there exists a function $u \geq 0$ in $\wop$ satisfying
    \begin{equation*}
        \sup_{t \geq 0} \Jl(tu) < C_{\infty}.
    \end{equation*}
In particular, this implies that $\Theta_{\lambda}^{-} < C_{\infty}$.
\end{lemma}
\begin{proof}
    Let $\bar\Lambda_0$ be as defined in \eqref{Lambda0bn}, ensuring that $C_{\infty} > 0$ for all $\lambda \in (0, \Lambda_0)$. Without loss of generality, we assume that $\inf_{B_{r_0}} a(x) = m_{a} > 0.$ Consider the function $\vep$ and get 
    \begin{equation*}
        \Jl(t\veps) \leq \frac{t^{p}}{p} \|\veps\|_{\wop}^{p} + \frac{t^{q}}{q} \|\veps\|_{\wsq}^{q} \leq C (t^{p} + t^{q}).
    \end{equation*}
    Consequently, there exists some $t_0 \in (0,1)$ such that
    \begin{equation*}
        \sup_{0 \leq t \leq t_0} \Jl(t\veps) < C_{\infty}.
    \end{equation*}
    Next, define the function
    \begin{equation*}
        h(t) = \frac{t^p}{p}\|\veps\|_{\wop}^{p} + \frac{t^q}{q} \|\veps\|_{\wsq}^{q} - \frac{t^{\pt}}{\pt}.
    \end{equation*}
    Noting that $h(0) = 0$, that $h(t) > 0$ for small $t$, and that $h(t) < 0$ for sufficiently large $t$, there exists $t_\varepsilon > 0$ such that
    \begin{equation*}
        \max_{t > 0} h(t) = h(t_\varepsilon),
    \end{equation*}
    where $t_\varepsilon$ is determined by solving $h'(t_\varepsilon) = 0$,
    \begin{equation*}
        t_\varepsilon^{\pt-q}=(t_\varepsilon^{p-q}\|\veps\|_{\wop}+\|\veps\|_{\wsq})\leq C(1+t_\varepsilon^{p-q}) .
    \end{equation*}
    This implies $0\leq  t_\varepsilon < t_1$  for some $t_1>0.$
    Using the inequality
    \begin{equation*}
        \int_{\Omega} a(x) |\veps|^{\delta} \geq m_{a} \int_{B_{r_0}} |\veps|^{\delta},
    \end{equation*}
    we obtain the estimate
    \begin{equation*}
    \begin{aligned}
            \sup_{t \geq t_0} \Jl(t\veps) &\leq \sup_{t > 0} h(t) - \lambda \frac{t_0^\delta m_a}{\delta} \int_{B_{r_0}} |\veps|^{\delta} \\
            &= \frac{t_\varepsilon^p}{p}\|\veps\|_{\wop}^{p}+ \frac{t_\varepsilon^q}{q} \|\veps\|_{\wsq}^{q}-\frac{t_\varepsilon^{\pt}}{\pt} -\lambda\frac{t_0^\delta m_a}{\delta}\int_{B_{r_0}} |\veps|^{\delta}\\
        & \leq \sup_{t\geq 0}\left( \frac{t^p}{p}\|\veps\|_{\wop}^{p} -\frac{t^{\pt}}{\pt}\right)+ \frac{t_1^q}{q} \|\veps\|_{\wsq}^{q}-\lambda\frac{t_0^\delta m_a}{\delta}\int_{B_{r_0}} |\veps|^{\delta}.
    \end{aligned}
    \end{equation*}
    Defining $g(t) = \frac{t^p}{p} \|\veps\|_{\wop}^{p} -  \frac{t^{\pt}}{\pt}$, we observe that $g$ attains its maximum at $\widetilde{t} = \left(\|\veps\|_{\wop}^{p}\right)^{\frac{1}{\pt - p}}$. Consequently,
    \begin{equation*}
        \sup_{t \geq 0} g(t) = g(\widetilde{t}) = \frac{\left(\|\veps\|_{\wop}^{p}\right)^{N/p}}{N}\leq \frac{1}{N}\left(S_{p}^{N/p} + O(\varepsilon^{\frac{N-p}{p-1}})\right).
    \end{equation*}
    This and the fact $q\leq p$ lead to the final estimate
    \begin{equation*}
    \begin{aligned}
        \sup_{t \geq t_{0}} \Jl(t\veps) &\leq \frac{1}{N}S_{p}^{N/p} + O(\varepsilon^{\frac{N-p}{p-1}}) + O(\varepsilon^\mpqcr)-\lambda\frac{t_0^\delta m_a}{\delta}\int_{B_{r_0}} |\veps|^{\delta}
        \\
        &\leq \frac{1}{N}S_{p}^{N/p}  + O(\varepsilon^\mpqcr)-\lambda\frac{t_0^\delta m_a}{\delta}\int_{B_{r_0}} |\veps|^{\delta}.
    \end{aligned}
    \end{equation*}
Given \(\|v_\varepsilon\|_{L^\delta}^{\delta} = O(\varepsilon^\alpha)\) and \(\varepsilon = \lambda^\beta\), the bound  $\sup_{t \geq t_{0}} \Jl(t\veps) < C_\infty$  
holds if  
\[
\lambda^{\beta\mpqcr} - \lambda^{1+\beta\alpha} \leq -C_\delta\lambda^\frac{p}{p-\delta}.
\]  
The previous condition holds for sufficiently small \(\lambda\), if \(\alpha < \mpqcr\) and  $\frac{1}{\mpqcr-\alpha} < \beta < \frac{\delta}{\alpha(p-\delta)}.$ Both conditions are satisfied if \(\alpha p < \delta \mpqcr\). Thus, for \(\lambda \in (0, \bar\Lambda_{00})\), a sufficient condition ensuring the bound $\sup_{t >0} \Jl(t\veps) < C_\infty$ is  
\[
\alpha < \frac{\delta\mpqcr}{p}.
\]  
This holds under either  
\[
\max\left\{\frac{Np}{\mpqcr+N-p},\pt(1-1/p)\right\}<\delta<q
\]  
or  
\[
\delta<\pt(1-1/p), \quad s<1-\frac{1}{q}\left(\frac{N-p}{p-1}-N\left(1-\frac{q}{p}\right)\right).
\]  
Applying the fibering map analysis on $\Jl$, we deduce the existence of \(\hat{t} > 0\) such that \(\hat{t} \veps \in \Nl^{-}\), concluding that \(\Theta_\lambda^- < C_{\infty}\).  
\end{proof}
Now, using similar arguments as in the previous section, we establish the existence of two nonnegative nontrival solutions.
\\[1mm]
\textit{\underline{\textbf{Sketch of Proof of Theorem \ref{bnthm}} }}
\\[3mm]
    By Proposition \ref{criticallevelbn}, there exists a minimizing sequence $\{u_k\}$ for $\Nl$ which is also  a $(PS)_{\Theta_{\lambda}}$ sequence for $\Jl$. Applying Lemma \ref{thetalambdaplusbn} and Lemma \ref{psrangebn}, we conclude that there exists $u_\lambda\in\wopo$ such that $u_k\to u_\lambda$ strongly in $\wopo$ for $\lambda\in(0,\bar\Lambda_{0})$. Now using arguments similar to Theorem \ref{exsonesolcri} we get the existence of one nonnegative solution.
    \\
    The results of Lemma \ref{projlembn} and Proposition \ref{criticallevelbn} remain valid even when $\Nl$ is replaced by $\Nl^-$. Consequently, we obtain a minimizing sequence $\{u_k\} \subset \Nl^-$ satisfying
    \begin{equation*}
        \Jl(u_k) = \Theta_\lambda^- + o_k(1), \quad \text{and} \quad \jl'(u_k) = o_k(1),
    \end{equation*}
    implying that $\{u_k\}$ forms a $(PS)_{\Theta_\lambda^-}$ sequence for $\Jl$.
    By applying Lemmas \ref{psrangebn} and \ref{talentilibn}, there exists a function $v_\lambda \in \wopo$ such that $u_k \to v_\lambda$ in $\wopo$. Theorem \ref{cdeltathmbn} further guarantees that $(\jl'(v_\lambda), v_\lambda) = 0$. Now we follow the arguments of Theorem \ref{multsolcri} to get the existence of second nonnegative nontrivial solution.\hfill\qed

\bibliographystyle{abbrv}
	\bibliography{refnew}

\begin{thebibliography}{10}

\bibitem{ABC94}
A.~Ambrosetti, H.~Brezis, and G.~Cerami.
\newblock Combined effects of concave and convex nonlinearities in some elliptic problems.
\newblock {\em J. Funct. Anal.}, 122(2):519--543, 1994.

\bibitem{AGP96}
A.~Ambrosetti, J.~Garcia~Azorero, and I.~Peral.
\newblock Multiplicity results for some nonlinear elliptic equations.
\newblock {\em J. Funct. Anal.}, 137(1):219--242, 1996.

\bibitem{AC23}
C.~A. Antonini and M.~Cozzi.
\newblock Global gradient regularity and a hopf lemma for quasilinear operators of mixed local-nonlocal type.
\newblock {\em arXiv preprint arXiv:2308.06075}, 2023.

\bibitem{BM19}
M.~Bhakta and D.~Mukherjee.
\newblock Multiplicity results for {$(p,q)$} fractional elliptic equations involving critical nonlinearities.
\newblock {\em Adv. Differential Equations}, 24(3-4):185--228, 2019.

\bibitem{BDVV22b}
S.~Biagi, S.~Dipierro, E.~Valdinoci, and E.~Vecchi.
\newblock A brezis-nirenberg type result for mixed local and nonlocal operators.
\newblock {\em arXiv preprint arXiv:2209.07502}, 2022.

\bibitem{BDVV22}
S.~Biagi, S.~Dipierro, E.~Valdinoci, and E.~Vecchi.
\newblock Mixed local and nonlocal elliptic operators: regularity and maximum principles.
\newblock {\em Comm. Partial Differential Equations}, 47(3):585--629, 2022.

\bibitem{BMV24}
S.~Biagi, D.~Mugnai, and E.~Vecchi.
\newblock A {B}rezis-{O}swald approach for mixed local and nonlocal operators.
\newblock {\em Commun. Contemp. Math.}, 26(2):Paper No. 2250057, 28, 2024.

\bibitem{BV24}
S.~Biagi and E.~Vecchi.
\newblock Multiplicity of positive solutions for mixed local-nonlocal singular critical problems.
\newblock {\em Calc. Var. Partial Differential Equations}, 63(9):Paper No. 221, 45, 2024.

\bibitem{BV24b}
S.~Biagi and E.~Vecchi.
\newblock On the existence of a second positive solution to mixed local-nonlocal concave-convex critical problems.
\newblock {\em arXiv preprint arXiv:2403.18424}, 2024.

\bibitem{BF22}
N.~Biswas and F.~Sk.
\newblock On generalized eigenvalue problems of fractional (p, q)-laplace operator with two parameters.
\newblock {\em Proceedings of the Royal Society of Edinburgh Section A: Mathematics}, pages 1--46, 2022.

\bibitem{BD13}
D.~Blazevski and D.~del Castillo-Negrete.
\newblock Local and nonlocal anisotropic transport in reversed shear magnetic fields: Shearless cantori and nondiffusive transport.
\newblock {\em Physical Review E—Statistical, Nonlinear, and Soft Matter Physics}, 87(6):063106, 2013.

\bibitem{BT15}
V.~Bobkov and M.~Tanaka.
\newblock On positive solutions for {$(p,q)$}-{L}aplace equations with two parameters.
\newblock {\em Calc. Var. Partial Differential Equations}, 54(3):3277--3301, 2015.

\bibitem{BCS13}
C.~Br\"andle, E.~Colorado, A.~de~Pablo, and U.~S\'anchez.
\newblock A concave-convex elliptic problem involving the fractional {L}aplacian.
\newblock {\em Proc. Roy. Soc. Edinburgh Sect. A}, 143(1):39--71, 2013.

\bibitem{BN83}
H.~Br\'ezis and L.~Nirenberg.
\newblock Positive solutions of nonlinear elliptic equations involving critical {S}obolev exponents.
\newblock {\em Comm. Pure Appl. Math.}, 36(4):437--477, 1983.

\bibitem{BZ03}
K.~J. Brown and Y.~Zhang.
\newblock The {N}ehari manifold for a semilinear elliptic equation with a sign-changing weight function.
\newblock {\em J. Differential Equations}, 193(2):481--499, 2003.

\bibitem{CD15}
W.~Chen and S.~Deng.
\newblock The {N}ehari manifold for nonlocal elliptic operators involving concave-convex nonlinearities.
\newblock {\em Z. Angew. Math. Phys.}, 66(4):1387--1400, 2015.

\bibitem{CMM18}
W.~Chen, S.~Mosconi, and M.~Squassina.
\newblock Nonlocal problems with critical {H}ardy nonlinearity.
\newblock {\em J. Funct. Anal.}, 275(11):3065--3114, 2018.

\bibitem{CKSR12}
Z.-Q. Chen, P.~Kim, R.~Song, and Z.~Vondravcek.
\newblock Boundary {H}arnack principle for {$\Delta+\Delta^{\alpha/2}$}.
\newblock {\em Trans. Amer. Math. Soc.}, 364(8):4169--4205, 2012.

\bibitem{CPS14}
E.~Colorado, A.~de~Pablo, and U.~S\'anchez.
\newblock Perturbations of a critical fractional equation.
\newblock {\em Pacific J. Math.}, 271(1):65--85, 2014.

\bibitem{CM94}
D.~G. Costa and C.~A. Magalh\~aes.
\newblock Variational elliptic problems which are nonquadratic at infinity.
\newblock {\em Nonlinear Anal.}, 23(11):1401--1412, 1994.

\bibitem{SCGG19}
E.~D. da~Silva, M.~L.~M. Carvalho, J.~V. Gon{\c{c}}alves, and C.~Goulart.
\newblock Critical quasilinear elliptic problems using concave-convex nonlinearities.
\newblock {\em Ann. Mat. Pura Appl. (4)}, 198(3):693--726, 2019.

\bibitem{DFV24}
J.~V. da~Silva, A.~Fiscella, and V.~A.~B. Viloria.
\newblock Mixed local-nonlocal quasilinear problems with critical nonlinearities.
\newblock {\em J. Differential Equations}, 408:494--536, 2024.

\bibitem{FGU09}
D.~G. de~Figueiredo, J.-P. Gossez, and P.~Ubilla.
\newblock Local ``superlinearity'' and ``sublinearity'' for the {$p$}-{L}aplacian.
\newblock {\em J. Funct. Anal.}, 257(3):721--752, 2009.

\bibitem{FM22}
C.~De~Filippis and G.~Mingione.
\newblock Gradient regularity in mixed local and nonlocal problems.
\newblock {\em Math. Ann.}, 388(1):261--328, 2024.

\bibitem{Pai11}
F.~O. de~Paiva.
\newblock Nonnegative solutions of elliptic problems with sublinear indefinite nonlinearity.
\newblock {\em J. Funct. Anal.}, 261(9):2569--2586, 2011.

\bibitem{DFR19}
L.~M. Del~Pezzo, R.~Ferreira, and J.~D. Rossi.
\newblock Eigenvalues for a combination between local and nonlocal {$p$}-{L}aplacians.
\newblock {\em Fract. Calc. Appl. Anal.}, 22(5):1414--1436, 2019.

\bibitem{DGJ24}
R.~Dhanya, J.~Giacomoni, and R.~Jana.
\newblock Interior and boundary regularity of mixed local nonlocal problem with singular data and its applications.
\newblock {\em arXiv preprint arXiv:2411.18505}, 2024.

\bibitem{DPV12}
E.~Di~Nezza, G.~Palatucci, and E.~Valdinoci.
\newblock Hitchhiker's guide to the fractional {S}obolev spaces.
\newblock {\em Bull. Sci. Math.}, 136(5):521--573, 2012.

\bibitem{DPV23}
S.~Dipierro, E.~Proietti~Lippi, and E.~Valdinoci.
\newblock ({N}on)local logistic equations with {N}eumann conditions.
\newblock {\em Ann. Inst. H. Poincar\'e{} C Anal. Non Lin\'eaire}, 40(5):1093--1166, 2023.

\bibitem{DV21}
S.~Dipierro and E.~Valdinoci.
\newblock Description of an ecological niche for a mixed local/nonlocal dispersal: an evolution equation and a new {N}eumann condition arising from the superposition of {B}rownian and {L}\'evy processes.
\newblock {\em Phys. A}, 575:Paper No. 126052, 20, 2021.

\bibitem{DP97}
P.~Dr{\'a}bek and S.~I. Pohozaev.
\newblock Positive solutions for the {$p$}-{L}aplacian: application of the fibering method.
\newblock {\em Proc. Roy. Soc. Edinburgh Sect. A}, 127(4):703--726, 1997.

\bibitem{FSV15}
A.~Fiscella, R.~Servadei, and E.~Valdinoci.
\newblock Density properties for fractional {S}obolev spaces.
\newblock {\em Ann. Acad. Sci. Fenn. Math.}, 40(1):235--253, 2015.

\bibitem{GK22}
P.~Garain and J.~Kinnunen.
\newblock On the regularity theory for mixed local and nonlocal quasilinear elliptic equations.
\newblock {\em Trans. Amer. Math. Soc.}, 375(8):5393--5423, 2022.

\bibitem{GL23a}
P.~Garain and E.~Lindgren.
\newblock Higher {H}\"older regularity for mixed local and nonlocal degenerate elliptic equations.
\newblock {\em Calc. Var. Partial Differential Equations}, 62(2):Paper No. 67, 36, 2023.

\bibitem{GA91}
J.~Garc\'ia~Azorero and I.~Peral~Alonso.
\newblock Multiplicity of solutions for elliptic problems with critical exponent or with a nonsymmetric term.
\newblock {\em Trans. Amer. Math. Soc.}, 323(2):877--895, 1991.

\bibitem{GA87}
J.~P. Garc\'ia~Azorero and I.~Peral~Alonso.
\newblock Existence and nonuniqueness for the {$p$}-{L}aplacian: nonlinear eigenvalues.
\newblock {\em Comm. Partial Differential Equations}, 12(12):1389--1430, 1987.

\bibitem{GGM22}
J.~Giacomoni, A.~Gouasmia, and A.~Mokrane.
\newblock Discrete {P}icone inequalities and applications to non local and non homogenenous operators.
\newblock {\em Rev. R. Acad. Cienc. Exactas F\'is. Nat. Ser. A Mat. RACSAM}, 116(3):Paper No. 100, 21, 2022.

\bibitem{GKS20}
D.~Goel, D.~Kumar, and K.~Sreenadh.
\newblock Regularity and multiplicity results for fractional {$(p,q)$}-{L}aplacian equations.
\newblock {\em Commun. Contemp. Math.}, 22(8):1950065, 37, 2020.

\bibitem{Hai03}
Y.~Haitao.
\newblock Multiplicity and asymptotic behavior of positive solutions for a singular semilinear elliptic problem.
\newblock {\em J. Differential Equations}, 189(2):487--512, 2003.

\bibitem{MMP19}
S.~A. Marano, G.~Marino, and N.~S. Papageorgiou.
\newblock On a {D}irichlet problem with {$(p,q)$}-{L}aplacian and parametric concave-convex nonlinearity.
\newblock {\em J. Math. Anal. Appl.}, 475(2):1093--1107, 2019.

\bibitem{SVWZ22}
X.~Su, E.~Valdinoci, Y.~Wei, and J.~Zhang.
\newblock Regularity results for solutions of mixed local and nonlocal elliptic equations.
\newblock {\em Math. Z.}, 302(3):1855--1878, 2022.

\bibitem{SW15}
Y.~Wei and X.~Su.
\newblock Multiplicity of solutions for non-local elliptic equations driven by the fractional {L}aplacian.
\newblock {\em Calc. Var. Partial Differential Equations}, 52(1-2):95--124, 2015.

\bibitem{Wil96}
M.~Willem.
\newblock {\em Minimax theorems}, volume~24 of {\em Progress in Nonlinear Differential Equations and their Applications}.
\newblock Birkh\"auser Boston, Inc., Boston, MA, 1996.

\end{thebibliography}

\end{document}